\documentstyle[twoside]{article}

\newtheorem{thm}{\indent{\sc Theorem}}[section]
\newtheorem{defn}[thm]{\indent{\sc Definition}} 
 
\newtheorem{por}[thm]{\indent{\sc Porism}}
\newtheorem{prop}[thm]{\indent{\sc Proposition}} 
\newtheorem{cor}[thm]{\indent{\sc Corollary}}
\newtheorem{ex}{\indent {\sc Example}}[section]
\newtheorem{prob}[thm]{\indent{\sc Open Problem}\index{Problem{,} Open}}

\newcommand{\ssect}{\subsection}

\newcommand{\eqref}[1]{equation~(\ref{#1})}
\newcommand{\eref}[1]{\eqref{#1}}

\newcommand{~}{\nolinebreak[3] }

\newcommand{\thinr}[1]{\makebox[0pt][r]{$\scriptstyle #1$}}
\newcommand{\thinl}[1]{\makebox[0pt][l]{$\scriptstyle #1$}}
\newcommand{\th}{\raisebox{0.6ex}{th}}
\newcommand{\del}{\nabla}
\newcommand{\seper}{\nopagebreak \begin{center}
	\underline{\hspace{2in}}
	\end{center}}

\newcommand{\bold}[1]{{\em (#1)}\index{#1}}

 \newcommand{\proof}[1]{\proo{#1}$\Box $}
\newcommand{\proo}[1]{{\em Proof: }#1}
\renewcommand{\Box}{ \hspace{1cm}{\rule{1.2ex}{2ex}}}

\newcommand{\rn}[1]{\left\lfloor #1 \right\rceil}

\newcommand{\D}{{{\rm D}}} 
 
\newcommand{\I}{{\rm I}}

\newcommand{\adj}[1]{\mbox {\rm adj}(#1)}

\newcommand{\romcoeff}[2]{\rn{{#1 \atop #2}}}

\newcommand{\action}[2]{\left\langle #2 \right\rangle 
 	\!\!_{_{\scriptstyle #1}}}

\newcommand{\SB}[1]{\mbox{{\bf \scriptsize #1}}}

\renewcommand{\bigotimes}{\otimes}
\renewcommand{\tilde}[1]{\widetilde{#1}}
  
\newcommand{\fig}[2]{\begin{table}[htbp]
	\caption{#1} 
	\begin{center}\fbox{ \scriptsize \begin{tabular}#2 \end{tabular}}
	\end{center}
\end{table}}

\newcommand{\figx}[2]{\begin{figure}[htbp] \caption{#1} 
	\begin{center} \scriptsize \mbox{}#2 \end{center} \end{figure}}

\begin{document}
\include{title}

\begin{abstract}
We take advantage of the
combinatorial  interpretations of many sequences of polynomials of binomial
type   to define a sequence 
of symmetric functions corresponding to each sequence of polynomials of
binomial type. We derive many of the results of Umbral Calculus in this context
including a Taylor's expansion and a
binomial identity for symmetric functions. Surprisingly, the delta 
operators for all the sequences of binomial type correspond to the same
operator  on symmetric functions.
\seper
{\bf Les suites de fonctions sym\'etriques de type binomial}

On s'appuie ici sur les interpr\'etations combinatoires de nombreuses
suites de polyn\^omes de type binomial pour d\'efinir une suite de
fonctions sym\'etriques associ\'ee \`a chque suite de polyn\^omes de
type binomial. On retrouve dans ce cadre, de nombreaux r\'esultats du
calcul ombral, en particulier   une version de la formule de Taylor et
la formule d'identit\'e du bin\^ome pour les fonctions sym\'etriques.
On s'aper\,coit que les op\'erateurs differentiels de degr\'e un pour
toutes les suite de polyn\^omes de type binomial correspondent \`a un
op\'erateur unique sur les fonction sym\'etriques. 
 \end{abstract}

\begin{center}
{\it Dedicated to the memory of\\
Rabbi Selig Starr}
\end{center}

\tableofcontents

\part{Linear Sequences}
\section{Introduction}

Although it is well known that
many sequence of polynomials of binomial type $p_{n}(x)$ enumerate the number
of functions from an $n$-element set to an $x$-element set enriched with a
certain type of structure on each block. There we show that if one allow
pseudospecies then every sequence of polynomials of binomial type is of
this form. Moreover, by counting the enriched functions more carefully we
define sequence of symmetric functions of binomial type $p_{n}({\bf y})$. Using
these methods we rederive many classical results of the theory of symmetric
functions, and of umbral calculus as well as a few new ones. We define
a shift operator and prove the accompanying binomial theorem, and then 
classify the set of operators invariant under it. The algebra of
shift-invariant operators turns out to naturally isomorphic to the dual Hopf
algebra of symmetric functions. Finally, we extend all these ideas to bases
$p_{\lambda }({\bf y})$ through the use of  {\em genera}---a generalization
of Joyal's species.

It is thought that the ideas in this paper would make an excellent
introduction to umbral calculus, species,  and symmetric functions for a
beginning graduate student in combinatorics.

\renewcommand{\backslash}{\setminus}
\newcommand{\id}[1]{{\rm Id}_{#1}({\bf y})}
\renewcommand{\Pi}{{\rm proj}}
\newcommand{\cent}[1]{\multicolumn{1}{|c}{#1}}

\ssect{Combinatorial Interpretations}\label{Interpret}
To devise an Umbral Calculus on symmetric functions we must first study the
combinatorial properties of symmetric function. The linear sequence of
symmetric functions associated to each sequence of polynomials of binomial type
is closely related to the sets of functions which these polynomials
enumerate. 

Classically, a sequence of polynomials $p_{n}(x)$ (with $\deg(p_{n}(x))=n$) is
said to be of {\em binomial type}\index{Binomial Type} if
\begin{equation}\label{bt}
p_{n}(x+a)=\sum _{k=0}^{n}{n\choose k }p_{k}(a)p_{n-k}(x).
\end{equation}
In this case, $q_{n}(x) = p_{n}(x)/n!$ is said to be a sequence of divided
powers since it obeys the identity
$$ q_{n}(x+a)=\sum _{k=0}^{n} q_{k}(a) q_{n-k}(x). $$
For example, 
\fig{Examples 1--3}{{lp{1.3in}|p{1.3in}|p{1.3in}|p{1.3in}}
 & \cent{Example Number} & \cent{1} & \cent{2}  & \cent{3}  \\ 
	\hline
A& Species & $\SB{Deg}: E \mapsto \{E \}$ 
& $ \SB{Inj}: E\mapsto 	\{\emptyset \}$ if $|E|=0$,
	$\SB{Inj}: E\mapsto E$ if $|E|=1$, and
	$\SB{Inj}: E\mapsto \emptyset$ if $|E|>1$.
& $\SB{Lin}: E\mapsto $ complete orderings of the set $E$\\
	\hline
B& Name of a Function Enriched by~A& Function & Injection & Disposition \\
	\hline
C& Sequence of Polynomials of Binomial Type Enumerated by~B
& $x^{n}$
& $(x)_{n}= x(x-1)\cdots (x-n+1)$
& $(x)^{n} = x(x+1)\cdots (x+n-1)$\\
	\hline
D&  Operator Associated with~C
& $\D $ & $\Delta $ & $\del $\\
	\hline 
E&  Operator Conjugate to~C
& $\D $ & $\log (\I +\D )$ & $-\log (\I -\D )$ \\
	\hline
F& Linear Sequence of Symmetric Functions of Binomial Type Enumerated by~B
& ${\rm Id}_{n}(\SB{y}) = (y_{1}+y_{2}+\cdots )^{n} $
& $(\SB{y})_{n}$
& $(\SB{y})^{n}$\\
	\hline
G& Linear Sequence of Symmetric Functions of Divided Powers Enumerated by~B
 & ${\rm Id}_{n}(\SB{y})/n! =$ $\sum _{\lambda \vdash
n}m_{\lambda }(\SB{y})/\lambda !$
 & $e_{n}(\SB{y}) = \sum _{\scriptstyle \mu \in
{\cal P}^{*} \atop \scriptstyle \mu \vdash n} \prod
_{k=1}^{n}y_{\mu _{k}}$
 & $h_{n}(\SB{y})=\sum _{\lambda \vdash n}m_{\lambda }(\SB{y})$ \\
	\hline
H& Exponential Generating Function for~F
& $\exp \left(t\sum _{y\in x}y\right)$ 
& $\prod _{y\in X}(1+yt)$ 
& $\prod_{y\in X} (1-yt)^{-1}$ }
\fig{Examples 4--6}{{l|p{1.3in}|p{1.3in}|p{1.3in}|p{1.3in}}
& \cent{Example Number} & \cent{4} & \cent{5}  & \cent{6}  \\ \hline 
A& Species& $\SB{F}= \exp (\SB{T}): E\mapsto $ rooted forests on the set $E$
& $\exp (\SB{Lin}):E\mapsto $ assemblies of linear orders on the set $E$
& $F_{1}=\exp ({\bf x}\SB{Deg}):E\mapsto $  rooted forests of trees of length
at most one on the set $E$ \\
	\hline
B& Name of a Function Enriched by~A&
Reluctant Function & Laguerre Function & Inverse-Abel Function \\
	\hline
C& Sequence of Polynomials of Binomial Type Enumerated by~B
& $A_{n}(x) = x(x+n)^{n-1}$
& $L_{n}(-x) = $ $\sum _{k=1}^{n}{n-1\choose k-1 }\frac{n!}{k!}x^{k}$
& $\mu _{n}(x) = \sum _{k=1}^{n}k^{n-k}{n\choose k}x^{k}$\\
	\hline
D&  Operator Associated with~C
& $\D E^{-1}$ & $\D /(\D -\I )$ &\\
	\hline 
E&  Operator Conjugate to~C
& & $\D /(\D -\I )$ & $\D E^{-1}$ \\
	\hline
F& Linear Sequence of Symmetric Functions of Binomial Type Enumerated by~B
& $A_{n}(\SB{y})= n! \sum _{\lambda \vdash n}\left(\prod _{i}\frac{\lambda
_{i}^{\lambda _{i}-2}}{(\lambda _{i}-1)!} \right)m_{\lambda }(\SB{y})$
& $L_{n}(\SB{y})$
& $\mu _{n}(\SB{y})$\\
	\hline
H& Exponential Generating Function for~F
&  & $\exp \left(\sum_{y\in X} \frac{yt}{1-yt}\right)$ 
& $\exp \left(\sum_{y\in X} yte^{yt} \right)$}
\begin{enumerate}
\item the powers $x^{n}$,
\item the lower factorial $(x)_{n}$,
\item the upper factorial $(x)^{n}$, 
\item the Abel polynomials $A_{n}(x)$,
\item the LaGuerre polynomials $L_{n}(-x)$, and
\item the inverse-Abel polynomials $\mu _{n}(x)$ 
\end{enumerate}
are all sequences  of binomial type. We show that
they all have similar combinatorial interpretations. In fact, we show that this
is in 
a sense typical of all sequences of binomial type. They all count the number of
functions from an $n$-element set $N$ to an $x$-element set $X$ which are
{\em enriched}\index{Enriched Functions} with a ``structure'' of some sort on
each {\em  fiber}\index{Fiber} where the 
fibers $f^{-1}(y)$ are the inverse images of elements of the range of the
function. 
\setcounter{figure}{-1}
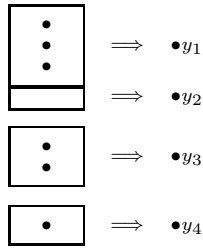
\figx{Typical Enriched Function\label{typ}}
{$\begin{array}{ccc}
 \framebox[1cm]{$\begin{array}{c}
\bullet \\ \bullet \\ \bullet 
\end{array} $} & \Longrightarrow & \bullet y_{1} \\[0.1in] 
\framebox[1cm]{$\begin{array}{c}
\end{array}$}&\Longrightarrow & \bullet y_{2}\\[0.1in]
 \framebox[1cm]{$\begin{array}{c}
 \bullet \\ \bullet 
\end{array} $} & \Longrightarrow & \bullet y_{3} \\[0.2in]	
\framebox[1cm]{$\begin{array}{c}
\bullet
\end{array}$} & \Longrightarrow & \bullet y_{4}
\end{array}$}
In Figure~\ref{typ} and each of the following figures, we display a typical
enriched function from the $n$-element set on the left to the $x$-element set
on the right. Note that we occasionally allow $x$ to be infinite.

Let us reformalize the above ideas in the language of species.
\begin{defn}[Species]
 Given a type of
structure (eg: rooted trees), its {\em species}\index{Species} ${\bf S}$ is the
functor from the category {\bf Sets} of finite sets and bijections to itself.
For any finite set $E$, we say that the members of the set ${\bf S}[E]$ are
${\bf S}$-structures, and for any bijection $f: E\rightarrow F$, we describe
the function  ${\bf S}[f]$ as a relabeling of ${\bf S}$-structures.
For this paper, we need to 
assume that there is only one structure on the empty set and there is at least
one structure on every one element set; that is,
$|{\bf S}[\emptyset]|=1$ and $|{\bf S}[\{0 \}]|\neq 0$.
\end{defn}

We define the sum of two species ${\bf S}_{1}$ and ${\bf S}_{2}$ on  $E$ by their disjoint
union 
$$ ({\bf S}_{1}+{\bf S}_{2})[E] = {\bf S}_{1}[E] \dot{\cup } {\bf S}_{2}[E].$$

Similarly, we define the product of ${\bf S}_{1}$ and ${\bf S}_{2}$ on $E$ to
be set of quadruples 
$$ ({\bf S}_{1}{\bf S}_{2})[E]=
 \left\{ 
\parbox{3.5in}{$(E_{1},E_{2} ,A_{1},A_{2})$: $E_{1}$ and $E_{2}$ are disjoint,
$E_{1}\cup E_{2}=E,$ and $A_{i}\in {\bf S}_{i}[E_{i}]$.}\right\}$$
In other words,   we divide $E$ in half,
 and place an ${\bf S}_{1}$-structure on one half, and an ${\bf
S}_{2}$-structure on the other.

Next, we define the exponentiation of a species ${\bf S}$ on $E$ to be
$$ \exp({\bf S})[E]=
 \left\{ 
\parbox{3.5in}{$(\phi , (A_{B})_{B\in \phi })$: $\phi$ is a partition of ${\bf
S}$, 
$(A_{B})_{B\in \theta }$ is a sequence of ${\bf S}$ structures one on each block of
$\phi $.}\right\}$$
That is, we divide $E$ into a number of parts and place an ${\bf S}$ structure on
each part.

Similarly, we define the composition of ${\bf S}_{1}$ with ${\bf S}_{2}$ on $E$ to be 
$$ \exp({\bf S})[\pi ]=
 \left\{ 
\parbox{3.5in}{$(\phi , (A_{B})_{B\in \phi },C)$: $\phi$ is a partition of
${\bf S}$, 
$(A_{B})_{B\in \theta }$ is a sequence of ${\bf S}$ structures one on each block of
$\phi $, and $C$ is an ${\bf S}$ structure on $\phi $.}\right\}$$
That is, we divide $E$ into a number of parts, place a ${\bf S}_{2}$-structure
on each parts, and finally place a ${\bf S}_{1}$-structure on the parts
themselves.  

Finally, we define the derivative ${\bf S}'$ of a species ${\bf S}$ on $E$ to
be the value of ${\bf S}$ on a set $E\dot{\cup} \{\infty  \}$ with one more
element. 

 Now, these enriched functions which we are counting are merely pairs
$(f,(a_{y})_{y\in X})$ where $f:N\rightarrow X$ and $a_{y}\in {\bf S}[f^{-1}(y)]$.

\figx{{Typical Function}\label{fun}}{$ \begin{array}{ccc}
\left. \begin{array}{c}
\bullet \\ \bullet \\ \bullet 
\end{array} \right\} 
& \Longrightarrow & \bullet y_{1} \\[0.1in]
& & \bullet y_{2}\\[0.1in]
\left.\begin{array}{c}
 \bullet \\ \bullet 
\end{array} \right\}
 & \Longrightarrow & \bullet y_{3} \\[0.1in]
{\bullet} & \Longrightarrow & \bullet y_{4}
\end{array} $}
{\sc Example \ref{Interpret}.1:} {\bf (Powers of $x$)} The powers $x^{n}$ are
the canonical example of a sequence 
of binomial type. They count all of the functions $f: N\rightarrow X$, since
each of the $n$ members of $N$ can be mapped independently onto any of the $x$
members of $X$.

These ordinary functions can be brought under the umbrella of the preceding
discussion by observing that they can be thought of as functions enriched by
the \index{Degenerate Species}{\em degenerate species:} ${\bf Deg}: E \mapsto
\{E \}$. Thus, on each fiber there is only one possible ``structure.''

\figx{{Typical Injection}\label{inj}}{$ \begin{array}{ccc}
{\bullet} & \Longrightarrow & \bullet y_{1}\\
&&\bullet y_{2}\\
{\bullet} & \Longrightarrow & \bullet y_{3}\\
{\bullet} & \Longrightarrow & \bullet y_{4}\\
&&\bullet y_{5}\\
&&\bullet y_{6}
\end{array} $}
{\sc Example \ref{Interpret}.2:} {\bf (Lower Factorial)}
The {\em lower factorial}\index{Lower Factorial} $(x)_{n}=x(x-1)\cdots (x-n+1)
=n!{x\choose n }$ counts
the number of injections of $N$ into $X$. Obviously, the first member of $N$
can be mapped in $x$ different way. This leaves $x-1$ choices for the second,
and so on.

Next, note that an injection is a function enriched with the species 
$$ {\bf Inj}: E\mapsto \left\{ \begin{array}{ll}
\{\emptyset \}&\mbox{if $|E|=0$,}\\[0.1in]
E&\mbox{if $|E|=1$, and }\\[0.1in]
\emptyset & \mbox{if $|E|>1$.}
\end{array} \right.$$
Thus, there are no structures available for ``illegal'' fibers containing two
or more points.

\figx{{Typical Disposition}\label{disp}\index{Disposition}}{$
\begin{array}{ccc} 
\framebox{\rule{0ex}{2ex}$\bullet < \bullet < \bullet $}
 & \Longrightarrow & \bullet y_{1} \\[0.1in]
\framebox{\rule{0ex}{2ex}}&\Longrightarrow & \bullet y_{2}\\[0.1in]
\fbox{\rule{0ex}{2ex}$ \bullet < \bullet $}
 & \Longrightarrow & \bullet y_{3} \\[0.1in]
\fbox{\rule{0ex}{2ex}$\bullet$} & \Longrightarrow & \bullet y_{4}
\end{array} $}
{\sc Example \ref{Interpret}.3:} {\bf (Upper Factorial)}
A {\em disposition}\index{Disposition} is a function enriched with the species {\bf Lin}
of linear orders; it is a function with a linear ordering on each of its
fibers. 

To count the number of dispositions, observe that if $N=\{1,2,3,\ldots ,n \}$,
then 1 can be mapped to any of the $x$ members of $X$. The same is true of 2;
however, if $f(1)=f(2)$, then we must also choose the order of 1 and 2. Hence,
there are total of $x+1$ choices. Regardless of which choice we take, there are
$x+2$ choices for 3, and so on. By a simple induction, the number of
dispositions is given by the {\em upper factorial}\index{Upper Factorial}
$(x)^{n}=x(x+1)(x+2)\cdots (x+n-1)$. 

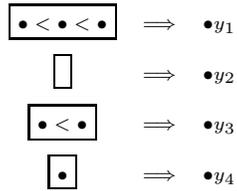
\figx{{Typical Reluctant Function}\label{rel}\index{Reluctant
Function}}{$ \begin{array}{ccc}
\framebox[1.5in]{\rule[-0.2in]{0cm}{0.4in}}
 & \Longrightarrow & \bullet y_{1} \\[0.15in]
\framebox[1.5in]{\rule[-0.15in]{0cm}{0.3in}}
 & \Longrightarrow & \bullet y_{2} \\[0.2in]	
\framebox[1.5in]{\rule[-0.3in]{0cm}{0.6in}}
 & \Longrightarrow & \bullet y_{3} \\[0.2in]	
\framebox[1.5in]{\rule[-0.25cm]{0cm}{0.5in}}
 & \Longrightarrow & \bullet y_{4} 
\end{array} $}
{\sc Example \ref{Interpret}.4:} {\bf (Abel Polynomials)}
We see that the {\em Abel polynomials}\index{Abel}
$A_{n}(x)=x(x+n)^{n-1}$ counts the number of functions enriched by the species
{\bf F} of labeled rooted forests. Such functions are called {\em reluctant
functions.}\index{Reluctant Functions}

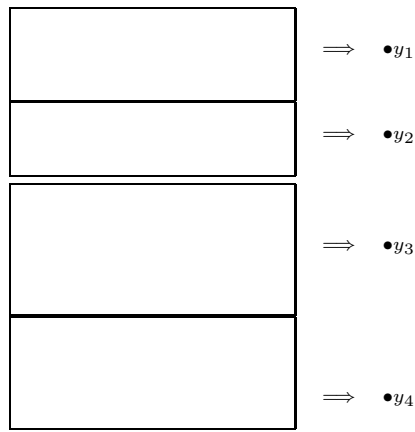
\figx{{Typical Laguerre Function}\label{lag}\index{Laguerre}}
{$ \begin{array}{ccc}
\framebox[1in]{$\begin{array}{c}
\bullet < \bullet < \bullet < \bullet\\
\bullet < \bullet \\
\bullet\\
\bullet
\end{array}$ }
 & \Longrightarrow & \bullet y_{1} \\[0.1in]
\framebox[1in]{$\bullet < \bullet < \bullet $}
 & \Longrightarrow & \bullet y_{2} \\[0.1in]
\framebox[1in]{}  & \Longrightarrow & \bullet y_{3} \\[0.1in]
\framebox[1in]{$\begin{array}{c}
\bullet < \bullet < \bullet \\
\bullet < \bullet < \bullet 
\end{array}$}
 & \Longrightarrow & \bullet y_{4} 
\end{array} $}
{\sc Example \ref{Interpret}.5:} {\bf (Laguerre Polynomials)}\index{Laguerre}
 Next, the {\em Laguerre polynomials} $$L_{n}(-x)= \sum _{k=1
}^{n} {n-1\choose k-1}\frac{n!}{k!}x^{k}$$ count the number of {\em
Laguerre functions}. These are functions from $N$ to $X$ enriched with the
species $\exp ({\bf Lin})$ of
collections of linear orders on each fiber. In a sense, the Laguerre functions
are related to the dispositions in the same way that functions are related to
injection. We will soon see the importance of this  relationship. 

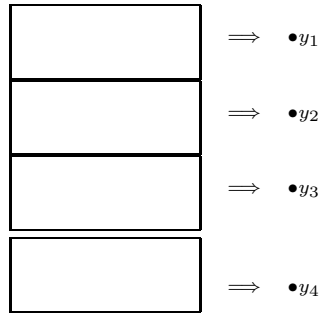
\figx{{Typical Inverse-Abel
Function}\label{ia}\index{Inverse-Abel Function} }
{$ \begin{array}{ccc}
\framebox[1in]{\rule[-0.15in]{0cm}{0.3in}}
 & \Longrightarrow & \bullet y_{1} \\[0.1in]
\framebox[1in]{\rule[-0.15in]{0cm}{0.3in}}
 & \Longrightarrow & \bullet y_{2} \\[0.1in]
\framebox[1in]{\rule[-0.15in]{0cm}{0.3in}}
 & \Longrightarrow & \bullet y_{3} \\[0.2in]	
\framebox[1in]{\rule[-0.15cm]{0cm}{0.3in}}
 & \Longrightarrow & \bullet y_{4} 
\end{array} $}
{\sc Example \ref{Interpret}.6:} {\bf (Inverse-Abel
Polynomials)}\index{Inverse-Abel Polynomial} 
Finally, we  observe that the {\em inverse-Abel
polynomials} $\mu_{n}(x)=\sum _{k=0}^{n}k^{n-k}{n\choose k}x^{k}$ counts the
number of functions enriched by the species $F_{1}$ of forests of rooted trees
of length at most one. 

\ssect{Symmetric Functions}\label{Functions}
To study these combinatorial relations more closely, we count 
enriched functions according to the size of the their fibers. 

Suppose we {\em  represent} an (enriched) function $f:N\rightarrow X$ by the
product  
$$ \tilde{f}=\prod _{i\in N}f(i). $$
For instance, the function depicted in Figure~\ref{typ} is represented by
$y_{1}^{3}y_{3}^{2}y_{4}.$ We {\em represent} a collection ${\cal F}$ of
enriched functions\index{Representation of a Function}
by the sum of the representations of the functions.
$$ \tilde{\cal F} =\sum _{f\in {\cal F}}\prod_{i\in N}f(i). $$
Thus, if ${\cal F}$ is invariant under permutation of $X$, then $\tilde{{\cal
F}}$ is a symmetric function. For example, if ${\bf S}$ is a species and ${\cal S}$
is the collection of all ${\bf S}$-enriched functions from $N$
 to $X=\{y_{1},y_{2},\ldots  \}$, then $\tilde{{\cal F}}$ is a linear sequence
symmetric functions
$p_{n}({\bf y)}$.

This linear sequence is said to be of {\em binomial type},\index{Binomial Type}
and its sister sequence $p_{n}({\bf y})/n!$ is said to be a linear sequence of
{\em divided powers}. To some extent, the sequence of divided powers enumerates
the set 
of functions from an ``unlabeled'' set $N$ to a labeled set $X$.
Both linear sequences are said to be {\em  associated}\index{Associated} with
the species ${\bf S}$, and {\em related}\index{Related} to the sequence of
polynomials  of binomial type $p_{n}(x)$ which enumerates the number of
${\bf S}$-enriched functions from an $n$-element set to an $x$-element set.

The homomorphism $\Pi $ makes this relationship quite explicit. It is defined
by setting the variables $y_{1},\ldots ,y_{x}$ equal to one, and setting all
other variables 
equal to zero.  Thus, $\Pi p_{n}({\bf y})=p_{n}(x)$

Note that $\Pi $ is characterized by its action on the monomial symmetric
functions 
$$ \Pi m_{\lambda }({\bf y})= \frac{(x)_{\ell (\lambda )}}{\prod _{i\geq
0}{\rm mult}_{i}{\lambda }!}$$
where
\begin{enumerate}
\item The {\em monomial symmetric function}\index{Monomial Symmetric
Function} is given by the sum
$$m_{\lambda}({\bf y})=\sum_{\alpha } y_{1}^{\alpha _{1}}
y_{2}^{\alpha _{2}}\cdots $$
over all distinct permutations $\alpha $ of the linear partition $\lambda$.
\item  A {\em linear partition}\index{Partition{, }Linear}
$\lambda $ is a nonincreasing infinite
sequence, $(\lambda  _{i})_{i\geq 1}$,  of nonnegative integers
which is eventually zero;
for example, $\lambda  =(17,2,2,1,0,0\ldots )$ is a linear partition.
A linear partition $\lambda $ is said to be a {\em partition of $n$}
if the sum of its parts is $n$, and we write $\lambda  \vdash n$.  
The set of all linear partitions is denoted by ${\cal P}$. 
Partitions may be compared in at least three ways.
\begin{enumerate}
\item For any partition $\lambda $ and vector $\alpha $,  if $\alpha _{i}\leq
\lambda _{i}$ for all $i$, then we write $\alpha  \leq 
\lambda $, and we denote by $\lambda -\alpha  $ their vector difference
$(\lambda -\alpha )_{i}=\lambda _{i}-\alpha _{i}$. Similarly, we denote by
$\lambda +\mu $ their vector sum.
\item Each nonzero $\lambda _{i}$ is called a {\em part}\index{Parts} $\lambda
$. In the 
above example,  the multiset of parts of $\lambda $ is $\{1,2,2,17\}$. The
number of parts of $\lambda $ is denoted $\ell (\lambda  )$, and the number of
parts of $\lambda $ equal to $i$ is denoted  ${\rm mult}_{i}(\lambda )$. If the
parts of $\mu $ form a submultiset of $\lambda $, that is to say if ${\rm
mult}_{i}(\mu )\leq {\rm mult}_{i}(\lambda )$ for all $i$; then we write $\mu
\subseteq \lambda $. We denote by $\lambda \backslash \mu $ the partition whose
multiset of parts is the difference between the multisets of parts for $\lambda
$ and $\mu $. Similarly, we denote by $\lambda \cup \mu $ the partition who
multiset of parts in the union of the multiset of parts for $\lambda $ and $\mu
$. 

A linear partition is said to have {\em distinct
parts}\index{Parts,Distinct} if 
its multiset of parts is, in fact, a set. The set of all linear
partitions with distinct parts is denoted by ${\cal P}^{*}.$
\item Finally, the Ferrers diagram of a partition $\lambda $  is defined to be
the set 
of ordered pairs $(i,j)$ such that $1\leq j\leq \lambda _{i}$. If the Ferrers
diagram 
of $\mu $ is a subset of the Ferrers diagram of $\lambda $, we write $\mu
\sqsubseteq \lambda $, and we denote by $\lambda /\mu $ the set difference
between the two Ferrers diagrams.
\end{enumerate}
\end{enumerate} 
$\Pi $ is obviously a homomorphism for the algebra of symmetric functions
$\Lambda $ to the algebra of polynomials ${\bf C}[x]$.

The monomial symmetric functions are the simplest known basis for $\Lambda $.
They allow us to give an explicit formula for $p_{n}({\bf y})$ in terms of its
species. 
\begin{prop}\label{explicit}
The linear sequence of divided powers associated
with the species ${\bf S}$ is given by the sum
$$ q_{n}({\bf y}) = \sum _{\lambda \vdash n}\left(\prod _{i} \frac{a_{\lambda
_{i}}}{\lambda _{i}!} \right) m_{\lambda }({\bf y}) $$
over partitions $\lambda $ of the nonnegative integer $n$
where  $a_{i}$ is the number of ${\bf S}$-structures on an $i$ element set.\Box
\end{prop}

The associated linear sequence of binomial type is thus
$$ p_{n}({\bf y}) = \sum _{\lambda \vdash n}{n\choose \lambda} \left(\prod _{i}
a_{\lambda _{i}}\right) m_{\lambda }({\bf y}) $$
where ${n\choose \lambda}=n!/\lambda _{1}!\lambda _{2}!\cdots $.
By assumption, $a_{0}=1$ and $a_{1}\geq 1$.

By projecting back to the polynomials, we derive a {\em new} result concerning
polynomials which gives the coefficients of one sequence of binomial type in
terms of the lower factorials.
\begin{cor}\label{lower pi}
Let $p_{n}(x)$ be the sequence of polynomials of binomial type associated with
the species $(a_{n})_{n\geq 1}$. Then
$$ p_{n}({\bf y}) = \sum _{\lambda \vdash n} {n\choose\lambda} \left(\prod
_{i}\frac{a_{\lambda_{i}}}{\mbox{mult}_{i}(\lambda )!} \right) (x)_{\ell
(\lambda )}.$$  
\end{cor}

\proof{This is a direct application of Proposition~\ref{explicit} bearing in
mind that the polynomial analog of the monomial symmetric function is $\Pi m_{\lambda
}({\bf y})= (x)_{\ell (\lambda )}/\prod_{i}\mbox{mult}_{i}(\lambda )!$.}

Let {\bf C} be the category of complex numbers along with ``maps'' from each
complex number $z$ to itself.
\figx{{Quasi-Species}\label{q species}}
{$ \begin{array}{ccc}
{\bf Sets} & {\SB{S}\atop \displaystyle \longrightarrow} & {\bf Sets}\\*
\thinr{\SB{Q}}\downarrow & \swarrow \thinl{\# }\\*
{\bf C}
\end{array} $}
Then the map $\#  : E \mapsto |E|$ is a functor from the category {\bf Sets} to
the category ${\bf C}$. We have seen that the linear sequence $p_{n}({\bf y})$
is completely determined by the composition of functors ${\bf Q}= \# \circ {\bf
S}$. We 
will therefore call any functor ${\bf Q}: {\bf Sets}\rightarrow $ {\bf C} a
{\em 
quasi-species} regardless of whether or not its splits into the composition of
$\#  $ and a species. In this case, $a_{n}$ above will refer to the value of
${\bf Q}$ on an $n$-element set. We define the linear sequence of 
divided powers associated with a quasi-species via Proposition~\ref{explicit}.
However, only in the case of species is there a clear combinatorial
interpretation as above.

NB: The examples below in this and the ensuing sections (through \S\label{bit})
are continuations 
of the examples in \S\ref{Interpret}. They are numbered according to their
section and their subject. For instance, Example~\ref{Functions}.1 is the
continuation of Example~\ref{Interpret}.1. 

{\sc Example \ref{Functions}.1:} {\bf (Powers Symmetric Function)}
The set of all functions from $N$ to $X$ is represented by
the linear sequence of symmetric function of binomial type
$$ \id{n}=(y_{1}+y_{2}+\cdots )^{n} =\sum _{\lambda \vdash
n}{n\choose \lambda  }m_{\lambda }({\bf y})$$
\fig{Powers Symmetric Function $\id{n}$}{{rcl}
${\rm Id}_{0}(a,b,c)$ &=& 1\\
${\rm Id}_{1}(a,b,c)$ &=& $a+b+c$\\
${\rm Id}_{2}(a,b,c)$ &=& $2ab+(a^{2}+b^{2})$\\
${\rm Id}_{3}(a,b,c)$ &=& $(a^{3}+b^{3}+c^{3})
+3(a^{2}b+a^{2}c+ba^{2}+bc^{2}+ca^{2}+cb^{2})+6abc$} 
whose sister sequence of divided powers is $(y_{1}+y_{2}+\cdots )^{n}/n!$. 
Thus, the symmetric functions ${\id{n}}$ {\em represents}  the
powers $x^{n}$. 

{\sc Example \ref{Functions}.2:} {\bf (Elementary Symmetric
Function)}\index{Lower Factorial}
The set of all injections from $N$ to $X$ is represented by
the lower factorial symmetric function
$$ ({\bf y})_{n}= n!e_{n}({\bf y}) $$
whose sister sequence of divided powers  $e_{n}({\bf y})$
otherwise known as the \index{Elementary Symmetric
Function}{\em elementary symmetric function} is (by Proposition~\ref{explicit})
given by the sum 
$$ e_{n}({\bf y}) = \sum _{\scriptstyle \mu \in
{\cal P}^{*} \atop \scriptstyle \mu \vdash n} \prod
_{k=1}^{n}y_{\mu _{k}}$$
over all linear partitions $\mu $ with
distinct parts.
\fig{Elementary Symmetric Function $e_{n}({\bf y})$}{{rcl}
$e_{0}(a,b,c)$ &=& 1\\
$e_{1}(a,b,c)$ &=& $a+b+c$\\
$e_{2}(a,b,c)$ &=& $ab+ac+bc$\\
$e_{3}(a,b,c)$ &=& $abc$}

{\sc Example \ref{Functions}.3:} {\bf (Complete Symmetric
Function)}\index{Upper Factorial} 
Dually, the set of all unlabeled dispositions\index{Dispositions} from $N$
to $X$ is represented by the {\em complete symmetric function}\index{Complete
Symmetric Function} $h_{n}(x)$.
\fig{Complete  Symmetric Function $h_{n}({\bf y})$}{{rcl}
$h_{0}(a,b,c)$ &=& 1\\
$h_{1}(a,b,c)$ &=& $a+b+c$\\
$h_{2}(a,b,c)$ &=& $ab+ac+bc +a^{2}+b^{2}+c^{2}$\\
$h_{3}(a,b,c)$ &=& $a^{3}+b^{3}+c^{3}
+a^{2}b+a^{2}c+ba^{2}+bc^{2}+ca^{2}+cb^{2} +abc$}   
 By Proposition~\ref{explicit}, the complete
symmetric function is given by the sum
$$ h_{n}({\bf y}) = \sum _{\scriptstyle \rho \in
{\cal P} \atop \scriptstyle \rho \vdash n} \prod
_{k=1}^{n}y_{\rho_{k}}$$
 over all linear partitions $\rho $. Similarly, the set of all dispositions is
represented by the upper factorial symmetric function $({\bf y})^{n}=n!h_{n}({\bf
y})$. 

{\sc Example \ref{Functions}.4:} {\bf (Abel Symmetric Function)}
The linear sequence of symmetric functions of binomial type for the set of Abel
functions from $N$ to $X$
 is given by the Abel symmetric function\index{Abel Symmetric Function} 
$$ A_{n}({\bf y})=n!\sum _{\lambda \vdash n}\left(\prod _{i}
\frac{\lambda_{i}^{\lambda_{i}-2}}{(\lambda_{i}-1)!} \right) 
m_{\lambda }({\bf y}).  $$  
\fig{Abel Symmetric Function $A_{n}({\bf y})$}{{rcl}
$A_{0}(a,b,c)$ &=& 1\\
$A_{1}(a,b,c)$ &=& $a+b+c$\\
$A_{2}(a,b,c)$ &=& $2(ab+ac+bc)+2(a^{2}+b^{2}+c^{2})$\\
$A_{3}(a,b,c)$ &=& $9(a^{3}+b^{3}+c^{3})
+6(a^{2}b+a^{2}c+ba^{2}+bc^{2}+ca^{2}+cb^{2}) +6abc$}   

\ssect{Generating Functions}\index{Generating Function}\label{rate}

We also wish to derive the exponential generation function for linear sequences
of symmetric functions of binomial type.
Suppose that ${\bf S}$ is a (quasi-)species such 
that ${\bf S}[E]$ has $a_{i}$ elements where $i$ is the number of elements in
the set $E$. By assumption, $a_{0}=1$ and $a_{1}\geq 1$.
The {\em generating function}\index{Generating Function} of the
species ${\bf S}$ is defined to be 
$$ {\rm Gen}[{\bf S}](t)=\sum _{i=0}^{\infty }a_{i} {t^{i}}/{i!}. $$
Notice that all the operations among species which we have defined are held
equally by their generating functions. That is,
\begin{eqnarray*}
\mbox{Gen}[{\bf S}_{1}](t)+\mbox{Gen}[{\bf S}_{2}](t)=\mbox{Gen}[{\bf
S}_{1}+{\bf S}_{2}](t)\\ 
\mbox{Gen}[{\bf S}_{1}](t)\mbox{Gen}[{\bf S}_{2}](t)=\mbox{Gen}[{\bf S}_{1}{\bf
S}_{2}](t)\\ 
\exp (\mbox{Gen}[{\bf S}_{2}](t))=\mbox{Gen}[\exp ({\bf S})](t)\\
\mbox{Gen}[{\bf S}_{1}](\mbox{Gen}[{\bf S}_{2}](t))=\mbox{Gen}[{\bf S}_{1}({\bf
S}_{2})](t) \\
d\mbox{Gen})[{\bf S}](t)/dt = \mbox{Gen}[{\bf S}'](t)
\end{eqnarray*}
The linear sequence of polynomials for a given species may be expressed in
terms of this generating function.
\begin{thm}
Let $p_{n}({\bf y})$ be the linear sequence of polynomials associated with the
quasi-species ${\bf S}=(a_{i})_{i\geq 1}$. Then $p_{n}({\bf y})$ is given
implicitly by the generating function
$$ \sum _{n=0}^{\infty } p_{n}({\bf y})t^{n}/n!=\prod _{y\in {\bf X}}{\rm Gen}[{\bf S}](yt).$$
\end{thm}

\proo{Consider the following sequence of equalities:
\begin{eqnarray}
 \sum _{n=0}^{\infty } p_{n}({\bf y})\frac{t^{n}}{n!} 
&=& \frac{1}{n!} \sum_{n=0}^{\infty } \label{a}
\sum _{\scriptstyle f:\{1,2,\ldots ,n \}\rightarrow X \atop \scriptstyle
\mbox{\tiny enriched by $\scriptstyle \SB{S}$}} \prod _{i=1}^{n} f(i)t\\*
&=& \frac{1}{n!} \sum _{n=0}^{\infty }
\sum _{\pi \vdash N} \label{b}
\sum _{\scriptstyle {\scriptstyle f:\{1,2,\ldots ,n \}\rightarrow X \atop
\scriptstyle \mbox{\tiny enriched by $\scriptstyle \SB{S}$}} \atop \scriptstyle
\mbox{\bf \tiny ker}(f)=\pi}  \prod _{y\in X} (yt)^{|\pi _{y}|} \\
\label{c} &=& \sum _{n=0}^{\infty }
\sum _{\pi \vdash N}
\prod _{y\in X} a_{|\pi _{y}|} (yt)^{|\pi _{y}|}/|\pi _{y}|!\\
\label{d} &=& \sum _{n=0}^{\infty }
\sum _{\alpha \vdash x}
\prod _{y\in X} a_{\alpha_{y}} (yt)^{\alpha_{y}}/\alpha _{y}!\\
\label{e} &=& \sum _{\alpha} \prod _{y\in X} a_{\alpha_{y}|}
(yt)^{\alpha_{y}} / \alpha_{y}!\\ 
 &=& \prod_{y\in X} \sum_{k=0}^{\infty } a_{k}(yt)^{k}/k!\nonumber \\*
&=& \prod_{y\in X} {\rm Gen}[\SB{S}](yt)\label{gen1}
\end{eqnarray}
where:
\begin{itemize}
\item in \eref{b} and \eref{c}, $\pi $ is a sequence of $x$  (possibly
empty) disjoint  sets indexed by $X$ whose union is $N$,
\item in \eref{a}, $f$ is a function from $\{1,2,\ldots ,n \}$ to $X$
enriched by the species ${\bf S}$; in \eref{b}, $f$ is in addition required
to have a fiber structure {\em (kernel)}\index{kernel} which corresponds to
$\pi $, and 
\item in \eref{d} and \eref{e}, $\alpha $ is a sequence of nonnegative
integers which total to $n$.\Box 
\end{itemize}}

Thus, by the remarks prior to Example \ref{Functions}.1, 
\begin{cor}
The sequence of polynomials of binomial type  associated
with the quasi-species ${\bf S}$ is given implicitly by the generating function
\begin{equation}\label{gen2}
\sum _{n=0}^{\infty } p_{n}(x)t^{n}/n! = {\rm Gen}[{\bf S}](t)^{x} = \exp (x\log {\rm
Gen}[{\bf S}](t)).\Box 
\end{equation}
\end{cor}

Note that ${\rm Gen}[{\bf S}](t)$ has a logarithm since $a_{0}=1$. Now, by a
classical  result of Umbral Calculus (see\cite{ch4}), $p_{n}(x)$ is the
conjugate sequence of   polynomials of binomial type for the delta operator
$\log ({\rm Gen}[{\bf S}]({\bf D} 
))$ where ${\bf D}$ is the derivative, and is associated and basic for the
compositional inverse of that delta operator.

Conversely, the sequence of polynomials of binomial type conjugate to the delta
operator $g({\rm D} )$ enumerates functions enriched by the quasi-species with
generating function $\exp(g(t)).$

Now, we see that counting enriching function enriched with {\em collections} of
${\bf S}$-structures as opposed to counting functions enriched with merely a
single 
${\bf S}$-structure is tantamount in terms of polynomials to Umbral composition with
the lower factorial sequence.\index{Lower Factorial}

In general, let $p_{n}({\bf y})$ and $q_{n}({\bf y})$ be the linear
sequences of symmetric functions which enumerate functions enriched with
the species ${\bf S}$ and $\exp (T)$  (collections of $T$-structures)
respectively. Next, let $p_{n}(x)=\Pi p_{n}({\bf y})$ and $q_{n}(x)=\Pi
q_{n}({\bf y})$ be their polynomial counterparts. Now, suppose $r_{n}(x)$ is 
the Umbral composition of $p_{n}(x)$ with $q_{n}(x)$, and that $r_{n}({\bf y})$
is the corresponding linear sequence of symmetric functions. Then $r_{n}({\bf
y})$ is associated with the species ${\bf S}({\bf T})$. That is, it enumerates functions
with a collection of ${\bf T}$-structures on each fiber and a single ${\bf
S}$-structure 
which unites the ${\bf S}$-structures. All of this will be better understood
later 
through the use of transfer operators.

{\sc Example \ref{rate}.1:} {\bf (Powers of $x$)}
The generating function for the degenerate species is 
$e^{t}$, and the logarithm of $e^{t}$ is $t$. Thus, the
generating\index{Generating Function} function for 
its associated linear sequence of symmetric functions of binomial type is 
$$ \exp \left(t\sum_{i\geq 0}y_{i} \right),$$
and the generating function for its sequence of polynomials of binomial type is
$ e^{xt}$. Thus, $x^{n}$ is the conjugate
graded sequence for ${\rm D}$. ${\rm D}$ is its own compositional inverse, so
$\D x^{n}=nx^{n} .$

{\sc Example \ref{rate}.2:} {\bf (Lower Factorial)}\index{Lower Factorial}
 The generating function for the species ${\bf Inj}$ is
$1+t$. The compositional inverse of $\log (1+t)$ is $e^{t}-1$, so the forward
difference \index{Forward Difference}
operator $e^{{\bf D}}-1$ is the delta operator for $(x)_{n}$. That is,
$(x+1)_{n}-(x)_{n}=n(x)_{n-1}$. 
By \eref{gen1}, the generating function\index{Generating Function} for the
elementary symmetric function\index{Elementary Symmetric Function} is given by 
$$ \prod _{n= 1}^{\infty }(1+y_{n}t)=\sum _{n=0}^{\infty }
e_{n}({\bf y})t^{n}. $$ 

{\sc Example \ref{rate}.3:} {\bf (Upper Factorial)}
The generating function for the species {\bf Lin} of linear orders is
$1/(1-t)$. The composition inverse of $\log (1/(1-t))$ is $1-e^{-t}$, so the
backward difference operator\index{Backward Difference Operator} $1-e^{{\bf
D}}$ is the delta operator for 
$(x)^{n}$. That is, $(x)^{n}-(x-1)^{n}=(x)^{n-1}.$
By \eref{gen1}, the generating function\index{Generating Function} for
the complete symmetric 
function\index{Complete Symmetric Function} is given by
$$ \prod _{n=1}^{\infty } (1-y_{n}t)^{-1}=\sum _{n=0}^{\infty }
h_{n}({\bf y})t^{n}. $$ 

{\sc Example \ref{rate}.4:} {\bf (Abel Function)}\index{Abel Function}
Let $T$ be the species of rooted trees. A tree with
its root removed is a rooted forest, and a collection of rooted trees is a
rooted forest. Thus, by the theory of species
$$ {\rm Gen}[T](x)=x{\rm Gen}[F](x), $$
and
$$ {\rm Gen}[F](x)=\exp ({\rm Gen}[T](x)). $$
Hence, the compositional inverse of the generating function for $F$ is
$x/e^{x}$. Now, $A_{n}(x)$ is the associated sequence of the Abel
operator\index{Abel Operator}
$\D e^{-\D }$, so as we claimed earlier $A_{n}(x)$ counts function enriched
with forests of rooted trees.

{\sc Example \ref{rate}.5:} {\bf (Laguerre Polynomials)}\index{Laguerre
Polynomials} The generating function for the species $\exp ({\bf Lin})$ of
collections of linear orders is $\exp (x/(1-x))$.
\fig{Laguerre Symmetric Function $L_{n}({\bf y})$}{{rcl}
$L_{0}(a,b,c)$ &=& 1\\
$L_{1}(a,b,c)$ &=& $a+b+c$\\
$L_{2}(a,b,c)$ &=& $3(ab+ac+bc)+2(a^{2}+b^{2}+c^{2})$\\
$L_{3}(a,b,c)$ &=& $7(a^{3}+b^{3}+c^{3})
+9(a^{2}b+a^{2}c+ba^{2}+bc^{2}+ca^{2}+cb^{2}) +6abc$}   
Thus, the Laguerre symmetric functions are given by the generating function 
$$\sum _{n=0}^{\infty } L_{n}({\bf y})\frac{t^{n}}{n!} = \prod_{y\in X}
\exp\left( \frac{yt}{1-yt}\right) = 
\exp\left( \sum_{y\in X} \frac{yt}{1-yt}\right).$$
Also, the Laguerre polynomials $L_{n}(-x)$ form the conjugate sequence for the
delta operator $\D /(1-\D )$. The compositional inverse of this operator is $\D
/(1+ \D )$. Thus, the Laguerre polynomials are associated and basic for the
 Weirstrass operator.\index{Weirstrass Operator}

{\sc Example \ref{rate}.6:} {\bf (Inverse-Abel
Polynomials)}\index{Inverse-Abel}  
The generating function for the species ${\bf T}_{1}$ of trees of
length at 
most one is $xe^{x}$ since there are exactly $n$ such trees on every set of
size $n$. Thus, the generating function for the species ${\bf F}_{1}=\exp ({\bf
T}_{1})$ of forests of such trees is $\exp (xe^{x})$.
\fig{Inverse-Abel Symmetric Function $\mu_{n}({\bf y})$}{{rcl}
$\mu _{0}(a,b,c)$ &=& 1\\
$\mu _{1}(a,b,c)$ &=& $a+b+c$\\
$\mu _{2}(a,b,c)$ &=& $2(ab+ac+bc)+2(a^{2}+b^{2}+c^{2})$\\
$\mu _{3}(a,b,c)$ &=& $3(a^{3}+b^{3}+c^{3})
+6(a^{2}b+a^{2}c+ba^{2}+bc^{2}+ca^{2}+cb^{2}) +6abc$}   
Hence, the generating function for the inverse-Abel symmetric functions is
$$\sum _{n=0}^{\infty } \mu_{n}({\bf y})\frac{t^{n}}{n!} = \prod_{y\in X} 
\exp\left( yte^{yt}\right)= \exp\left(\sum _{y\in X}yte^{yt}\right).$$

\ssect{The Symmetric Derivative}\index{Derivative,Symmetric}\index{Symmetric
Derivative}\label{derives}
We define a derivation on the ring of symmetric functions. 
(Up to a constant) this operator 
fills the role of the basic delta operator  simultaneously for each
linear sequence of symmetric functions associated with a 
sequence of polynomials of binomial type.

We define the {\em symmetric derivative} as follows.
\begin{equation}\label{sd}
 \D p(y_{1},y_{2},\ldots ) = \lim_{\epsilon \rightarrow 0} \frac{p(\epsilon
,y_{1},y_{2},\ldots )-p(0,y_{1},y_{2},\ldots )}{\epsilon } 
= \frac{{\bf d}p(\epsilon ,y_{1},y_{2},\ldots )}{{\bf d}\epsilon }. 
\end{equation}
In other words, we ``add'' a new variable, take the derivative with respect to
it, and set it equal to zero.

For calculations involving $\D $, it is useful to note how $\D $ behaves with
respect to the monomial symmetric functions. 
$$ \D m_{\lambda }({\bf y}) =\left\{ \begin{array}{ll}
m_{\lambda \backslash 1}({\bf y}) & \mbox{if 1 is a part of $\lambda$, and }\\*
0&\mbox{otherwise.}
\end{array} \right. $$

{\sc Example \ref{derives}.1:} {\bf (Powers Symmetric Function)}
The symmetric derivative acts on ${\rm Id}_{n}({\bf y})$ as the derivative acts
on the powers of $x$. $$\D {\rm Id}_{n}({\bf y})=n{\rm Id}_{n-1}({\bf y}).$$ 

{\sc Example \ref{derives}.2:} {\bf (Lower Factorial Symmetric Function)} Also,
$$\D e_{n}({\bf y})=e_{n-1}({\bf y}),$$
so  $$\D ({\bf y})_{n}=n({\bf y})_{n-1}.$$
 Hence, in this case the symmetric derivative plays
the role of the forward difference operator\index{Forward Difference Operator}
rather than that of the ordinary derivative. 
(We could also have derived this using the fact $e_{n}({\bf
y})=m_{(1^{n})}({\bf y}).$)

{\sc Example \ref{derives}.3:} {\bf (Upper Factorial Symmetric Function)}\index{Upper
Factorial} 
Also, 
\begin{equation}\label{eD}
\D h_{n}({\bf y})=h_{n-1}({\bf y}),
\end{equation}
 so $\D
({\bf y})^{n}=n({\bf y})^{n-1}$. Again, $\D $ play
the role of the backwards difference operator.

We can use \eref{eD} to calculate the symmetric derivative of the Schur
function. The {\em Schur function}\index{Schur function} is given by the
determinant
\begin{equation}\label{Schur Definition}
s_{\lambda }({\bf y})={\rm det}\left[h_{\lambda _{i}-i+j}({\bf y})
\right]_{i,j= 1}^{\infty }.
\end{equation}
Now, since $\D $ is a derivation, 
\begin{eqnarray*}
\D s_{\lambda }({\bf y}) &=& \sum _{n= 1}^{\infty }{\rm
det}\left[h_{\lambda_{i}-i+j -\delta _{in}} ({\bf y}) \right]_{i,j=
0}^{\infty}\\*  
&=& \sum _{n\in \lambda } s_{\lambda\backslash n \cup (n-1)}({\bf y}) 
\end{eqnarray*}
where the last sum is over the {\em distinct} parts of the partion $\lambda $.
For example, 
$$ \D s_{5221}({\bf y})=s_{4221}({\bf y})+s_{5211}({\bf y}) + s_{522}({\bf
y})$$

By analogous reasoning, for $\mu \sqsubseteq \lambda $, the {\em skew Schur
function}\index{Schur Function, 
Skew}\index{Skew Schur Function}
$$ s_{\lambda /\mu } = {\rm det} \left[h_{\lambda _{i}-\mu_{j}-i+j}\right]
({\bf y})$$ 
has a derivative of
$$ \D s_{\lambda /\mu }({\bf y}) = \sum _{n\in \lambda } s_{(\lambda\backslash
n \cup (n-1))/\mu }({\bf y}).$$

\seper 
Considering the combinatorial interpretation, we see that the only contribution
arises from functions
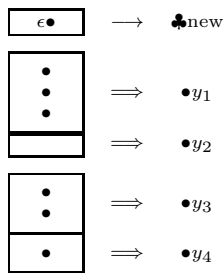
\figx{{Typical Contribution to the Derivative}\label{derfig}}
{$ \begin{array}{ccc}
\framebox[1cm]{$\epsilon \bullet$}& \longrightarrow & \clubsuit
\mbox{new}\\[0.1in] 
 \framebox[1cm]{$\begin{array}{c}
\bullet \\ \bullet \\ \bullet 
\end{array} $} & \Longrightarrow & \bullet y_{1} \\[0.1in]
\framebox[1cm]{$\begin{array}{c}
\end{array}$}&\Longrightarrow & \bullet y_{2}\\[0.1in]
 \framebox[1cm]{$\begin{array}{c}
 \bullet \\ \bullet 
\end{array} $} & \Longrightarrow & \bullet y_{3} \\[0.1in]
\framebox[1cm]{$\begin{array}{c}
\bullet
\end{array}$} & \Longrightarrow & \bullet y_{4}
\end{array} $}
 in which the new variable has a fiber of size one. This
fiber can contain any of the $n$ elements of $N$, and any of the $a_{1}$ many
structures allowed for one-element sets. Thus,
$$\D p_{n}({\bf y})=na_{1}p_{n-1}({\bf y}),$$
or equivalently, 
$$ \epsilon (\D ^{m}p_{n})=n!a_{1}^{n}\delta _{nm} $$
where $\epsilon $
acts on a  symmetric function $p({\bf y})$ by mapping all of the variables to
zero: $\epsilon p({\bf y})=p(0,0,\ldots )$; $\epsilon $ is the symmetric analog
of evaluation at zero since 
\begin{equation}\label{zero}
\epsilon p({\bf y})=[\Pi p({\bf y})]_{x=0}.
\end{equation}

Hence, $a_{1}^{-1}\D$ plays the role of the basic delta operator. This is true
even in the case of linear sequences of symmetric functions arising from
quasi-species. 

{\sc Example \ref{derives}.4:} {\bf (Abel Symmetric Function)}\index{Abel
Symmetric Function} 
$\D A_{n}({\bf y})=nA_{n-1}({\bf y}).$

{\sc Example \ref{derives}.5:} {\bf (Laguerre Symmetric
Function)}\index{Laguerre Symmetric Function} 
$\D L_{n}({\bf y})=nL_{n-1}({\bf y}).$

{\sc Example \ref{derives}.6:} {\bf (Inverse-Abel Symmetric
Function)}\index{Inverse-Abel Symmetric Function} 
 $\D \mu _{n}({\bf y})=n\mu _{n-1}({\bf y}).$

\ssect{The Iterated Symmetric Derivative}\label{mail}

Note further that if we define the {\em symmetric $i$-th
derivative}\index{Derivative,Symmetric}\index{Symmetric $i\th$ derivative}
$$ \D_{i} m_{\lambda }({\bf y}) =\left\{ \begin{array}{ll}
i!m_{\lambda \backslash i}({\bf y})&\mbox{if i is a part of $\lambda$, and }\\*
0&\mbox{otherwise,}
\end{array} \right. $$
then 
$$\D_{i} p_{n}({\bf y})=(n)_{i}a_{i}p_{n-i}({\bf y})$$
for any linear sequence of symmetric functions of binomial type $p_{n}({\bf
y})$. In particular,
$$ \D _{i}p_{n}({\bf y})=a_{i}{i!}\D ^{i}p_{n}({\bf y})/a_{1}^{i}. $$
Nevertheless, $\D_{i}$ is not a multiple of $\D ^{i}$.

More intuitively, $\D _{i}p({\bf y})$ is calculated by introducing a new
variable $y_{0}$, and then differentiating $p_{y_{0},y_{1},y_{2}\ldots }$ $i$
times with respect to $y_{0}$ and setting $y_{0}$ to zero.

Except when $i=1$, $\D _{i}$ is not a derivation. For example, 
$$\D_{2} m_{(1)}({\bf y})^{2} = \D _{2}\left( 2m_{(11)}({\bf y})
+m_{(2)}{({\bf y})} \right)= 1$$
whereas 
$$ 2m_{(1)}({\bf y})\left(\D _{2}m_{(1)}({\bf y}) \right)=0. $$

Note that $\D _{i}$ is $\D _{s_{i}({\bf y})}$ in the notation of
A.~Lascoux~\cite{L}.\index{Lascoux{,} Alain}

{\sc Example \ref{mail}.1:} {\bf (Powers Symmetric Function)} This is the
action of the 
Iterated derivative on the Powers Symmetric Function. $\D _{i}
\id{n}= (n)_{i}\id{n-i}.$

{\sc Example \ref{mail}.2:} {\bf (Lower Factorial Symmetric Function)}
\index{Lower Factorial} 
We would also like to compute the action of the Iterated derivative on the
elementary and complete symmetric functions. \index{Elementary Symmetric
Function} 
\begin{eqnarray*}
\D _{i}e_{n}({\bf y}) &=& \D _{i} ({\bf y})_{n}/n!\\*
&=& (n)_{i} \delta _{i,1} ({\bf y})_{n-i}/n!\\*
&=& \delta _{i,1} e_{n-i}({\bf y}).
\end{eqnarray*} 
Thus, $\D _{i}e_{n}({\bf y})=0$ except in the case of the ordinary symmetric
derivative where $\D _{1}e_{n}({\bf y})=e_{n-1}({\bf y}).$

{\sc Example \ref{mail}.3:} {\bf (Upper Factorial Symmetric
Function)}\index{Upper Factorial} Similarly, \index{Complete Symmetric
Function} 
\begin{eqnarray*}
\D _{i}h_{n}({\bf y}) &=& \D _{i} ({\bf y})^{n}/n!\\*
&=& (n)_{i} i! ({\bf y})^{n-i}/n!\\*
&=& i! h_{n-i}({\bf y}).
\end{eqnarray*} 

\ssect{The Binomial Theorem}\label{bit}

We next wish to devise an analog of the usual binomial identity for linear
sequences of polynomials of binomial type (\eref{bt}). However, we must first
determine the appropriate analog of the shift operator $E^{a}$. 
We define the {\em symmetric shift operator}\index{Symmetric Shift
Operator} $E^{a}$ by 
$$E^{a}=\sum _{n= 0}^{\infty }a^{n}\D _{n}/n!.$$ 
Hence,
$$ E^{a}m_{\lambda }({\bf y}) = \sum_{n=0}^{\infty } a^{n}m_{\lambda /n}({\bf
y})=m_{\lambda }(a,y_{1},y_{2},\ldots ), $$ 
and thus
\begin{equation}\label{shifty}
E^{a}p({\bf y}) =  p(a,y_{1},y_{2},\ldots ).
\end{equation}

Considering \eref{sd}, this is a very appropriate definition of the
symmetric shift operator. Note that by \eref{shifty}, $E^{a}$ is a
isomorphism of the field of symmetric functions.
However, note that its inverse is not $E^{-a}$. In fact, $E^{a}E^{b}\neq
E^{a+b}$ except in the case where $a$ or $b$ equals zero, for
$$E^{a}E^{b}m_{(11)}({\bf y}) = m_{(11)}({\bf y}) + (a+b)m_{(1)}({\bf y}) +
ab$$
whereas 
$$E^{a+b}m_{(11)}({\bf y}) = m_{(11)}({\bf y}) + (a+b)m_{(1)}({\bf y}).$$

Also, in consideration of \eref{sd}, the operators $E^{a}$ all
commute. 
For any set $S$, we define the shift $E^{S}$ to be the composition in any order
of the 
various shifts $E^{a}$ for all $a\in S$. In the case of an infinite set $S$,
this definition makes sense via the use of inverse limits. Thus, for ${\bf
z}=\{z_{1},z_{2},\ldots  \}$,
$$ E^{\SB{z}}p({\bf y})=p({\bf z}\cup {\bf y}). $$

Returning to our theory of linear sequence of symmetric functions of binomial
type. We can now derive the binomial identity.
\begin{thm}[Binomial Theorem]\index{Binomial Theorem}
If $p_{n}({\bf y})$ is a linear sequence of symmetric functions of binomial
type, then for all complex numbers $a$
\begin{equation}\label{bineq}
p_{n}(a,y_{1},y_{2},\ldots )=\sum _{k=0}^{n} {n\choose k}p_{k}(a,0,0,\ldots
) p_{n-k}({\bf y}). 
\end{equation} 
\end{thm}

{\bf Proof:} The left side of \eref{bineq} enumerates the set of
enriched functions from $N$ 
to $X\cup \{a \}$. The right side of \eref{bineq}
 counts the number of ways to choose a subset
of $N$
\figx{{Binomial Theorem}\label{binfig}\index{Binomial Theorem}}
{$ \begin{array}{ccc}
 \framebox[1cm]{$\begin{array}{c}
\spadesuit \\ \spadesuit \\ \spadesuit
\end{array} $} & \longrightarrow & \clubsuit a \\[0.1in]
\framebox[1cm]{$\begin{array}{c}
\end{array}$}&\Longrightarrow & \bullet y_{1}\\[0.1in]
 \framebox[1cm]{$\begin{array}{c}
 \bullet \\ \bullet 
\end{array} $} & \Longrightarrow & \bullet y_{2} \\[0.1in]
\framebox[1cm]{$\begin{array}{c}
\bullet
\end{array}$} & \Longrightarrow & \bullet y_{3}
\end{array} $}
 (represented in Figure~\ref{binfig} by spades $\spadesuit $) and map it to
$a$ with an enriched function (represented by a thin arrow
$\longrightarrow$) while mapping the 
remainder of $N$ (represented by dots $\bullet $) into $X$ with another
enriched function (represented by a thick arrow $\Longrightarrow$). Obviously,
these two sets of enriched functions are identical.$\Box $ 

\begin{cor}\label{repeated}\index{Roman Sequence}
If $p_{n}({\bf y})$ is a linear sequence of symmetric functions of binomial
type, and ${\bf y}$ and ${\bf z}$ are sets of variables, then
$$ p_{n}({\bf y\cup {\bf z}})= \sum _{k=0}^{n} {n\choose k} p_{n-k}({\bf
y})p_{k}({\bf z}). \Box $$
\end{cor}

\proof{By iterating Theorem 2 $m$-times, we have 
\begin{equation}\label{starn}
p_{n}(z_{1},\ldots ,z_{m},{\bf y})=\sum _{k}{n \choose k} p_{n-k}({\bf y})
\sum_{k_{1},\ldots ,k_{m}} {k \choose k_{1},\ldots ,k_{m} } \prod _{i=1}^{m} 
p_{k_{i}}(z_{i},0,0,\ldots ).
\end{equation}
where the inner sum is over sequences of nonnegative integers summing to $k$.
However, this sum may be computed  by \eref{starn} itself (with $m$
having a value of one less).
$$\sum _{k_{1},\ldots
,k_{m}} {k \choose k_{1},\ldots ,k_{m} } \prod _{i=1}^{m}
p_{k_{i}}(z_{i},0,0,\ldots )=p_{k}(z_{1},\ldots ,z_{m},0,0,\ldots )
$$   
Thus, 
$$ p_{n}(z_{1},\ldots ,z_{m},{\bf y})=\sum _{k=0}^{n}{n \choose k}
p_{k}(z_{1},\ldots ,z_{m},0,0,\ldots )p_{n-k}({\bf y}).
 $$
Our result now follows for infinitely many variables by the usual technique of
inverse limits.}

{\sc Example \ref{bit}.1:} {\bf (Powers Symmetric Function)}
 Applying the symmetric version of the binomial theorem to
$\id{n}$ we rediscover the usual binomial theorem
$$ (a+y_{1}+y_{2}+\cdots )^{n}=\sum _{k=0}^{n}{n\choose
k}a^{k}(y_{1}+y_{2}+\cdots )^{n-k}.$$

{\sc Example \ref{bit}.2:} {\bf (Lower Factorial Symmetric
Function)}\index{Lower Factorial} 
 Applying the binomial theorem to the lower factorial
symmetric function, we derive the following identity held be the elementary
symmetric function:\index{Elementary Symmetric Function}
$$ e_{n}(a,y_{1},y_{2},\ldots ) = e_{n}(y_{1},y_{2},\ldots ) +
ae_{n-1}(y_{1},y_{2},\ldots ), $$
or more generally
$$ e_{n}(z_{1},z_{2},\ldots ,y_{1},y_{2},\ldots ) = \sum _{k=0}^{n}
e_{k}(y_{1},y_{2},\ldots ) e_{n-k}(z_{1},z_{2},\ldots ). $$

{\sc Example \ref{bit}.3:} {\bf (Upper Factorial Symmetric
Function)}\index{Upper Factorial} 
 Applying the binomial theorem to the upper factorial
symmetric function, we derive the following identity held be the complete
symmetric function:\index{Complete Symmetric Function}
$$ h_{n}(a,y_{1},y_{2},\ldots ) = \sum _{k=0}^{\infty }
a^{k}h_{n-k}(y_{1},y_{2},\ldots ).$$

\seper 

We are also interested in the shifts of the other important symmetric
functions. 
The {\em power sum symmetric function}\index{Power Sum Symmetric Function}
$\mbox{pow}_{n}({\bf y})$ is defined by  
\begin{equation}\label{power}
\mbox{pow}_{n}({\bf y})=m_{(n)}({\bf y})=\sum _{i\geq 1}y_{i}^{n}.
\end{equation}
Clearly, $\mbox{pow}_{n}({\bf y}\cup {\bf z})=\mbox{pow}_{n}({\bf y})+
\mbox{pow}_{n}({\bf y})$.
Note that the kernel of the homomorphism $\Pi $ from the symmetric functions to
the polynomials is described very simply in terms of the power sum symmetric
function. Since $\Pi \mbox{pow}_{n}({\bf y})=x$, the kernel of $\Pi $ is the
algebra generated by the symmetric functions $\mbox{pow}_{n}({\bf
y})-\mbox{pow}_{m}({\bf y})$ where $m$ and $n$ are distinct integers.
 
The shift of the monomial symmetric function is given by
$$ m_{\lambda }({\bf y}\cup {\bf z}) = \sum _{\mu \cup \nu =\lambda } m_{\mu
}({\bf y}) m_{\nu }({\bf z}). $$ 

Finally, we can compute the shift of the Schur function. The Schur function can
be defined by \eref{Schur Definition} or by the ratio of anti-symmetric
functions in $n$ variables 
$$s_{\lambda }({\bf y})= a_{\lambda +\delta }({\bf y})/a_{\delta }({\bf y})$$
where $a_{\mu}({\bf 
y}) = {\rm det}\left[y_{i}^{\lambda _{j}+n-j}\right]$ and $\delta
=(n,n-1,\ldots ,2,1)$. 
This definition is consistent with changes of variables, since the ratio is
commutes with the map taking 
one of the variables to zero. After taking inverse limits, this definition is
equivalent to \eref{Schur Definition}. Hence, 
$a_{\delta }({\bf y})$ is simply the Vandermonde\index{Vandermonde}
determinant. 
$$ a_{\delta }({\bf y})=\prod_{1\leq i<j\leq n}(y_{i}-y_{j}). $$
Thus,
$$ E^{a}a_{\delta }({\bf y}) = (a-y_{1})(a-y_{2})\cdots (a-y_{n})a_{\delta
}({\bf y}).$$
Expanding by minors, we have 
\begin{eqnarray*}
E^{a}s_{\lambda }({\bf y})&=& {\rm det}\left[\begin{array}{cccc}
a^{\lambda _{1}+n+1}&a^{\lambda _{2}+n} & \cdots  & a_{\lambda _{n}+1}\\
y_{1}^{\lambda _{1}+n+1} & y_{1}^{\lambda _{2}+n} & \cdots  &
y_{1}^{\lambda_{n}+1}\\  
y_{2}^{\lambda _{1}+n+1} & y_{2}^{\lambda _{2}+n} & \cdots  &
y_{2}^{\lambda_{n}+1}\\  
\vdots & \vdots & \ddots & \vdots \\
y_{n}^{\lambda _{1}+n+1} & y_{n}^{\lambda _{2}+n} & \cdots  &
y_{n}^{\lambda_{n}+1}
\end{array}
\right]\\
&=&\sum _{j=1}^{n+1} (-1)^{j}a^{n-j+\lambda _{j}}({\bf y}) a_{\delta +(\lambda
_{1},\lambda _{2},\ldots, \lambda_{i-1}, \lambda
_{i+1}+1,\lambda_{i+2}+1,\ldots ,\lambda _{n+1}+1)}({\bf y})\\ 
&=& \left(\prod _{i=1}^{n+1} \frac{a}{a-y_{i}} \right) \left(\sum _{i=1}^{n+1}
a^{\lambda _{i}-i} s_{\lambda
_{1},\lambda _{2},\ldots, \lambda_{i-1}, \lambda
_{i+1}+1,\lambda_{i+2}+1,\ldots ,\lambda _{n+1}+1}({\bf y}) \right).\\
\sum _{\mu }a^{|\lambda |-|\mu |} s_{\mu }({\bf y})
\end{eqnarray*}
where the sum is over partitions $\mu \sqsubseteq \lambda $ such that 
the difference  of there Ferrers diagrams\index{Tableaux}\index{Ferrers
Diagram} $\lambda -\mu $ is a collection of disjoint horizontal strips.

This fact leads immediately to the definition of the Schur function in terms of
tableaux \cite{Mac}, and thence to the identities
\begin{eqnarray*}
s_{\lambda }({\bf y}\cup {\bf z}) &=& \sum _{\mu\sqsubseteq \lambda  }
s_{\lambda /\mu }({\bf y}) s_{\mu }({\bf y})\\
s_{\lambda/\nu}({\bf y}\cup {\bf z}) &=& \sum _{\mu \sqsubseteq \lambda }
s_{\lambda /\mu }({\bf y}) s_{\mu /\nu }({\bf y}).	
\end{eqnarray*}
Also, we see that the polynomial analog of the Schur function $\Pi s_{\lambda
}({\bf y})$ is the number of standard tableaux of shape $\lambda $ utilizing
the alphabet $1,2,\ldots ,x$.

There is more work yet to be done regarding the shift operator. For example,
\begin{prob}
What is the relation between  the action of the shift operator $E^{a}$ over the
ring of symmetric function and  the theory of Baxter algebras?
\end{prob}

\ssect{Shift-Invariant Operators}\index{Shift-Invariant}\label{shift}

Notice that  the symmetric derivatives commute with each other and with the
symmetric shift operator. These are essentially the only linear operators which
do so.

\begin{thm}[Classification of Shift-Invariant Operators] 
Let $\theta $ be a linear operator on symmetric functions. $E^{a}\theta =\theta
E^{a}$ for all $a$ if and only if $\theta $ is a complex formal power series in
the iterated derivatives $\D_{i}$.
\end{thm}

\proo{{\bf (If)} Immediate.

{\bf (Only If)} Let $b_{\lambda \mu }$ be the coefficients of the action of
$\theta $ on the monomial symmetric functions
$$ \theta m_{\lambda }({\bf y})=\sum _{\mu \in {\cal P}} b_{\lambda \mu }
m_{\mu }({\bf y}).  $$
Then 
\begin{eqnarray*}
\theta E^{a}m_{\lambda }({\bf y}) &=& \theta \sum _{i=0}^{\infty }
a^{i}m_{\lambda \backslash i}({\bf y})\\
&=& \sum _{\mu \in {\cal P}}\sum _{i=0}^{\infty } a^{i}b_{\lambda \backslash i,
\mu }m_{\mu }({\bf y})\\ 
E^{a}\theta m_{\lambda }({\bf y}) &=& E^{a} \sum _{\mu\in {\cal P} } b_{\lambda
\mu }m_{\mu }({\bf y}) \\
&=& \sum _{\mu \in {\cal P}}\sum _{i=0}^{\infty} a^{i}b_{\lambda \mu }m_{\mu
\backslash i}({\bf y}). 
\end{eqnarray*}
We equate coefficients of $a^{i}m_{\mu }({\bf y})$, and determine that
$$ b_{\lambda \backslash i,\mu }=b_{\lambda ,\mu\cup (i)}. $$
Thus, $b_{\lambda \mu }$ must equal zero unless  $\mu \subseteq \lambda $.
Moreover, $b_{\lambda \mu }$ depends only on  $\lambda \backslash
\mu $.  

Hence, the spaces of linear shift invariant operators is contained in the span
of the operators 
$$ \theta _{\nu }:m_{\lambda }({\bf y})\mapsto m_{\lambda \backslash \nu }({\bf
y}). $$ 
Thus, all shift-invariant operators are formal
power series in terms of the symmetric derivatives, since
$$\theta_{\nu }= \prod_{i=0}^{\infty } {\D_{\nu_{i}}}/{\nu _{i}!}.\Box $$  }

\begin{por}
The iterated derivatives $\D _{i}$ are algebraically independent.$\Box $
\end{por}

Let us write $\D _{\lambda }$ for the product of derivatives $\theta _{\lambda
}= \prod_{i}\D_{\lambda _{i}}/\lambda _{i}!$. The $\D _{\lambda }$ form a basis for the
module of linear shift-invariant operators. 
Clearly, we have the following identity
\begin{equation}\label{cast}
\epsilon \D _{\lambda }m_{\nu }({\bf y}) = \delta _{\lambda \nu }.
\end{equation}
Thus, the operator $\epsilon \D _{\lambda }$ is equal to the
classical operator $\langle h_{\lambda }({\bf y})|$ where $\langle | \rangle$
is the inner product defined by 
\begin{equation}\label{star}
\action{}{h_{\lambda }({\bf y})|m_{\lambda }({\bf y})}=\delta _{\lambda\mu}.
\end{equation}

\begin{thm}[Expansion Theorem] \index{Expansion Theorem} \label{sex}
Let $\theta $ be a linear shift invariant operator. Then
$$ \theta = \sum _{\lambda \in {\cal P}} 
\left( \epsilon\theta m_{\lambda }({\bf y})\right) \D _{\lambda }  .\Box $$ 
\end{thm}

\begin{thm}[Taylor's Theorem] \index{Taylor's Theorem}\label{stall}
Let $p({\bf y})$ be a symmetric function. Then
$$ p({\bf y}) = \sum _{\lambda \in {\cal P}} 
\left( \epsilon \D_{\lambda }p({\bf y}) \right) m_{\lambda }({\bf y}).\Box $$ 
\end{thm}

Also, we have the symmetric analog of Roman's identity for formal power series
of binomial type.

\begin{thm}[Roman's Identity]\index{Roman Sequence}
Let $p_{n}({\bf y})$ be a linear sequence of symmetric functions of binomial
type, and let $\theta $ and $\phi$ be shift invariant linear operators. Then
\begin{equation}\label{romeq}
\epsilon \theta \phi p_{n}({\bf y})= \sum_{k=0}^{n} {n\choose k} (\epsilon
\theta p_{k}({\bf y}))(\epsilon \phi p_{n-k}({\bf y})).
\end{equation}
\end{thm}

\proof{It  suffices to consider the case $\theta =\D
_{1}^{b_{1}}\D 
_{2}^{b_{2}}\cdots $ and $ \phi =\D _{1}^{c_{1}}\D _{2}^{c_{2}}\cdots .$
However, then both sides of \eref{romeq} equal zero unless 
$$ n=(b_{1}+c_{1})+2(b_{2}+c_{2})+\cdots  $$
in which case both sides equal
$$ n!a_{1}^{b_{1}+c_{1}}a_{2}^{b_{2}+c_{2}} $$
where $(a_{i})_{i\geq 0}$ is the quasi-species associated with $p_{n}({\bf
y}).$}

Thus, we have the following generalization of the symmetric version of binomial
theorem. 

\begin{cor}\label{forgotten}\index{Binomial Theorem}
Let $p_{n}({\bf y})$ be a linear sequence of symmetric functions of binomial
type, and let $\theta $ be a shift invariant linear operator. Then
$$ \theta p_{n}({\bf y}) = \sum_{k=0}^{n} {n\choose k} (\epsilon \theta
p_{k}({\bf y}))p_{n-k}({\bf y}).\Box  $$ 
\end{cor}

\ssect{Coalgebras}

The foregoing theory may be profitably be recast in terms of coalgebras. Let
$\Lambda _{\SB{y}}$ represent the algebra of symmetric functions in the
variables ${\bf y}$. $\Lambda _{\SB{y}}$ is isomorphic for any choice of
variables ${\bf y}$, so we  write $\Lambda $ for the algebra of symmetric
functions in any set of variables.

Now, $\Lambda _{\SB{x}}\otimes \Lambda _{\SB{y}}$ is isomorphic to the
algebra of formal power series invariant under permutations of ${\bf x}$ or
${\bf y}$.  Thus, $\Lambda_{\SB{x}\cup \SB{y}}$ is naturally a subalgebra of
$\Lambda _{\SB{x}}\otimes \Lambda _{\SB{y}}$. Hence, $E=E^{\SB{y}}$ may be
thought of as a map from $\Lambda $ to $\Lambda \otimes \Lambda $, and
$\epsilon=\epsilon _{\SB{y}}$ is a map from $\Lambda _{\SB{y}}=\Lambda $ to the
complex numbers.

From this point of view, we arrive at a result due to L.
Geissinger\index{Geissinger{, }L.} \cite{G}.
\begin{thm}
The symmetric functions $\Lambda $ form a commutative Hopf algebra
when they are equipp\-ed with 
\begin{itemize}
\item  the shift operator $E: \Lambda \rightarrow \Lambda \bigotimes  \Lambda $ as
the comultiplication, 
\item  the evaluation operator $\epsilon : \Lambda  \rightarrow {\bf C}$ as the
counit,
\item  the usual multiplication $\mu : \Lambda \bigotimes \Lambda \rightarrow
\Lambda   $,
\item the inclusion 
$\iota :{\bf C}\rightarrow \Lambda $ as the unitary map, and 
\item the classical involution
$$ \begin{array}{rrcl}
\omega: & \Lambda &\rightarrow & \Lambda \\*
&h_{n}({\bf y})&\mapsto &e_{n}({\bf y}).
\end{array}$$
as the antipode. 
\end{itemize}
\end{thm}

\proo{}{\bf (Commutative Algebra)} $\Lambda $ is well know to be an algebra.
That is, it obeys the commutative diagrams
\figx{{Associativity}\label{associative}}
{$ \begin{array}{ccc}
\Lambda \bigotimes \Lambda & {\scriptstyle I\otimes\mu \atop \displaystyle
\longleftarrow} 
&\Lambda \bigotimes \Lambda \bigotimes \Lambda  \\*
\thinr{\scriptstyle \mu } \downarrow & & \downarrow\thinl{\scriptstyle \mu
\otimes I} \\* 
\Lambda &{\scriptstyle \mu \atop \displaystyle \longleftarrow} & \Lambda
\bigotimes \Lambda  
\end{array} $}
 in Figures \ref{associative}--\ref{commutative}
\figx{{Unitary Property}\label{unit}}
{$ \begin{array}{ccccc}
&&\Lambda \bigotimes \Lambda \\*
&\thinr{\iota \otimes I}\nearrow & \downarrow\thinl{\mu} &
\nwarrow\thinl{I\otimes\iota } \\*
{\bf C }\bigotimes \Lambda &\sim & \Lambda & \sim  & \Lambda \bigotimes {\bf C}
\end{array} $}
where
$\tau $ is the commutation map $p({\bf x})\bigotimes q({\bf x})\mapsto q({\bf
x}) \bigotimes p({\bf x})$. 
\figx{{Commutativity}\label{commutative}}
{$ \begin{array}{ccc}
\Lambda \bigotimes \Lambda  \\
\thinr{\tau }\updownarrow & \thinr{\mu }\searrow\\
\Lambda \bigotimes \Lambda &{\displaystyle \longrightarrow \atop \scriptstyle
\mu }  & \Lambda 
\end{array} $}

{\bf (Cocommutative Coalgebra)} Cocommutivity follows from commutivity once
we observe that $\Lambda $ is a bialgebra; however, for now we observe this
directly from Figure~\ref{cocommutative} since $p({\bf x}\cup {\bf y})= p({\bf
y}\cup {\bf x})$.
\figx{{Cocommutivity}\label{cocommutative}}
{$ \begin{array}{ccc}
\Lambda_{\SB{y}} \bigotimes \Lambda_{\SB{x}} & {\tau \atop \displaystyle
\longleftrightarrow} &\Lambda_{\SB{x}} \bigotimes\Lambda_{\SB{y}}\\* 
\thinr{\scriptstyle E^{\SB{x}} } \uparrow & & \uparrow\thinl{\scriptstyle
E^{\SB{y}}}\\*  
\Lambda_{\SB{y}} & \sim & \Lambda_{\SB{x}}
\end{array} $}
To actually, show that $\Lambda $ is a coalgebra refer to Figures
\ref{coassociative} and \ref{counit}. 
\figx{{Coassociativity}\label{coassociative}}
{$ \begin{array}{ccc}
\Lambda_{\SB{x}} \bigotimes \Lambda_{\SB{y}} & {\scriptstyle I\otimes
E^{\SB{z}} \atop \displaystyle \longrightarrow} 
&\Lambda_{\SB{x}} \bigotimes \Lambda_{\SB{y}} \bigotimes \Lambda_{\SB{z}}  \\*
\thinr{\scriptstyle E^{\SB{y}} } \uparrow & & \uparrow\thinl{\scriptstyle
E^{\SB{y}}\otimes I} \\* 
\Lambda_{\SB{x}} &{\scriptstyle E^{\SB{z}} \atop \displaystyle \longrightarrow}
& 
\Lambda_{\SB{x}} \bigotimes \Lambda_{\SB{z}}
\end{array} $}
Coassociativity merely means that $p({\bf
x}\cup ({\bf y}\cup {\bf z}))=p(({\bf x}\cup {\bf y})\cup {\bf z})$, and the
counitary property means that $p(y_{1},y_{2},\ldots ,0,0,\ldots
)=p(y_{1},y_{2},\ldots )$.
\figx{{Counitary Property}\label{counit}}
{$ \begin{array}{ccccc}
&&\Lambda _{\SB{x}}\bigotimes \Lambda _{\SB{y}}\\*
&\thinr{\epsilon _{\SB{x}}\otimes I}\swarrow & \uparrow\thinl{E^{\SB{y}}} &
\searrow\thinl{I\otimes \epsilon _{\SB{y}}} \\*
{\bf C }\bigotimes \Lambda _{\SB{y}} &\sim & \Lambda_{\SB{x}} & \sim  &
\Lambda_{{\bf x}}\bigotimes {\bf C}  
\end{array} $}

{\bf (Bialgebra)} There are four requirements one of which is represented by
Figure~\ref{bialgebra}.
\figx{Bialgebra\label{bialgebra}}{$ \begin{array}{ccc}
\Lambda_{\SB{x}} \bigotimes \Lambda_{\SB{x}} & {E^{\SB{y}}\otimes E^{\SB{y}}
\atop \displaystyle 	
\longrightarrow} &\Lambda_{\SB{x}} \bigotimes \Lambda _{\SB{y}} \bigotimes
\Lambda_{\SB{x}}
\bigotimes \Lambda _{\SB{y}}\\* 
\thinr{\scriptstyle \mu _{\SB{x}} } \downarrow & &
\downarrow\thinl{\scriptstyle \mu_{\SB{x}}\otimes \mu _{\SB{y}}}\\*
\Lambda_{\SB{x}} & {\displaystyle \longrightarrow \atop \scriptstyle
E^{\SB{y}}} &  
\Lambda_{\SB{x}}\bigotimes \Lambda _{\SB{y}}
\end{array} $}
They are all satisfied since as we have already observed $E^{\SB{y}}$ and
$\epsilon _{\SB{y}}$ are both algebra homomorphism. 

{\bf (Hopf Algebra)} Finally, we must show that 
$$ \theta \omega_{\SB{x}}E^{\SB{y}}p({\bf x})=
\theta \omega_{\SB{y}} E^{\SB{y}}p({\bf x})= \epsilon _{\SB{x}}(p({\bf x}))$$
where $\theta $ is the evaluation of ${\bf x}$ at ${\bf y}$ and $\omega
_{\SB{y}}$ is the classical involution of $\Lambda _{\SB{y}}$.
However, this follows immediately from the following series of equalities
\begin{eqnarray*}
1 &=& \left(\prod _{y\in X} (1-yt) \right)\left(\prod _{y\in X} (1-yt)^{-1}
\right) \\
&=& \left(\sum _{n\geq 0} (-1)^{n}e_{n}({\bf y})t^{n} \right) \left(\sum_{n\geq
0} h_{n}({\bf y}) \right)\\ 
&=& \sum _{n\geq 0} \sum _{i=0}^{n} (-1)^{i}e_{i}({\bf y})h_{n-i}({\bf
y})t^{n}\\ 
 \delta _{n,0} &=&   \sum _{i=0}^{n} (-1)^{i}e_{i}({\bf y})h_{n-i}({\bf
y})t^{n}\\ 
&=&  \sum _{i=0}^{n} (-1)^{i}(\theta\omega h_{i}({\bf x}))h_{n-i}({\bf
y})t^{n}\\ 
&=& \theta \omega _{\SB{y}}E^{\SB{y}}h_{n}({\bf x}).\Box 
\end{eqnarray*}

In this notation, an operator $\theta $ is
shift-invariant\index{Shift-Invariant Operator} if and only if
$(\theta \otimes I)E^{\SB{x}} = E^{\SB{x}}\theta $

The dual Hopf algebra $\Lambda^{*}$ is generated by the maps $m_{\lambda
}^{*}: m_{\mu }({\bf y})\rightarrow \delta _{\lambda \mu}$. Thus, by
\eref{cast} $m_{\lambda}^{*} = \epsilon \D _{\lambda}$. Therefore,
automorphisms of the algebra of shift-invariant operators correspond to
automorphisms of the coalgebra of symmetric functions. These coalgebra maps
are called Transfer operators.\index{Transfer Operator}\index{Coalgebra
Map} 

It is best to compare the study of symmetric functions by coalgebra methods, to
the study of polynomials also by coalgebra methods. 
In that case, we are dealing with the usual multiplication along with the
comultiplication map
$$\Delta x^{n}=\sum _{k=0}^{n}{n\choose k}x^{k}\bigotimes x^{n-k},$$
the counitary map of evaluation at zero, and the antipodal map of substitution
by $-x$.
However, this Hopf algebra is seen to be a homomorphic image of the Hopf
algebra of symmetric functions as follows:
Define $\Pi _{\SB{x}}$ to act on the variables ${\bf x}$ and set exactly $x$
of them equal to 1, and the remainder to zero (where $x$ is a new
variable). As noted before $\Pi _{\SB{x}}$ is an algebra homomorphism from the
symmetric functions to the polynomials. Moreover, it is easy to see that 
$\Pi $ satisfies all of the necessary relations---for example it satisfies
$$\Delta \Pi _{\SB{x}}=(\Pi_{\SB{x}}\otimes \Pi _{\SB{y}} )\Delta $$
and \eref{zero}---so that $\Pi $ is a Hopf algebra homomorphism.

However, although the umbral calculus of polynomials and of symmetric functions
is best thought of in terms of coalgebras, we have no coalgebra version of the
logarithmic algebra.
\begin{prob}
Is the logarithmic algebra of \cite{ch4} naturally a coalgebra? 
\end{prob}
If so, then perhaps in analogy to polynomials we have
$\mu (X^{n}\bigotimes
X^{m}) =\romcoeff{n+m}{n}X^{n+m}$.
Would this then to be the incidence coalgebra for the poset of hybrid sets?
\cite{hybrid}

\part{Full Sequences}
\ssect{Introduction}

Since we require $\deg(p_{n}(x))=n$, every sequence of polynomials
$p_{n}(x)$ is a basis for the ring of polynomials. Clearly that is not the case
for linear sequences of symmetric functions. Since linear sequences of
symmetric functions are 
indexed by a single integer, there are insufficiently many. By the fundamental
theorem of symmetric functions, the linear sequences of elementary and complete
symmetric symmetric functions are both
transcendence bases. However, this is not true in general since the sequence 
$\id{n}=(y_{1}+y_{2}+\cdots )^{n}$ used in Example \ref{Functions}.1
is not a transcendence basis.\index{Transcendence Basis}

\ssect{Full Sequences}

Thus, we are led to discuss larger sequences of symmetric functions which are
indexed by partitions and span the space of symmetric functions.

\begin{defn}\bold{Full Sequence}\label{full}
Order the partitions in reverse lexicographical order as in \cite{Mac}. Then
define the {\em exact degree}\index{Degree} of a nonzero homogeneous symmetric
function $p({\bf y})$ of degree $n$ to be the largest partition $\lambda \vdash
n$ such that the coefficient $c_{\lambda }=\epsilon \Delta _{\lambda }p({\bf
y})$ is nonzero in the Taylor expansion (Theorem~\ref{stall})
$$ p({\bf y})=\sum_{\lambda \in {\cal P}} c_{\lambda }m_{\lambda }({\bf y}). $$

A {\em full sequence}\index{Full Sequence} then is a sequence $(p_{\lambda
}({\bf y}))_{\lambda \in {\cal P}}$ of homogeneous symmetric functions indexed
by partitions $\lambda $ such that each $p_{\lambda }({\bf y})$ is exactly of
degree $\lambda $.
\end{defn}

When $p_{\lambda}({\bf y})$ is a full sequence and $\alpha $ is a vector
of nonnegative integers with finite support, then we  denote by
$p_{\alpha}({\bf y})$ the symmetric function $p_{\mu }({\bf y})$ where $\mu $
is the partition formed by arranging the members of $\alpha $ in weakly
descending order. Moreover, we will adopt this convention for any sequence
indexed by partitions. 

Clearly, every full sequence is a basis for the algebra of symmetric functions.
Moreover,  the popular bases mentioned above---Schur $s_{\lambda }({\bf
y})$ and monomial $m_{\lambda }({\bf y})$---are full sequences.\cite{Mac}
Moreover, the linear sequences of complete $h_{n }({\bf y})$, elementary
$e_{n}({\bf y})$, and power $\mbox{pow}_{n}({\bf y})$ symmetric functions can
be extended to full sequences via the product rules:
\begin{eqnarray*}
e_{\lambda }({\bf y})&=&\prod _{i}e_{\lambda_{i}}({\bf y})\\*
h_{\lambda }({\bf y})&=&\prod _{i}h_{\lambda_{i}}({\bf y})\\*
\mbox{pow}_{\lambda }({\bf y}) &=&\prod_{i}\mbox{pow}_{\lambda _{i}}({\bf y})
\end{eqnarray*}
By the fundamental theorem of symmetric functions, these are bases for the
module of symmetric functions. In fact, it can be seen that $h_{\lambda }({\bf
y})$ and $e_{\lambda }({\bf y})$ are full sequences.

To determine the analog of $e_{\lambda }({\bf y})$, $h_{\lambda }({\bf y})$,
and $\mbox{pow}_{\lambda }({\bf y})$
in the theory of polynomials, we apply the homomorphism $\Pi$.
\begin{eqnarray*}
\Pi e_{\lambda }({\bf y}) &=& \prod _{(i,j)\in \lambda } (x-i+1)/i \\*
\Pi h_{\lambda }({\bf y}) &=& \prod _{(i,j)\in \lambda } (x+i+1)/i \\*
\Pi \mbox{pow}_{\lambda }({\bf y}) &=& x^{\ell (\lambda )}.
\end{eqnarray*}

Now, since $E^{a}$ is an isomorphism, we can also calculate the shift of these
full sequences. 
\begin{equation}\label{e shift}
 e_{\lambda }({\bf y}\cup {\bf z} ) = \sum _{\alpha}
e_{\alpha }({\bf y}) e_{\lambda -\alpha }({\bf z})
\end{equation}
where the sum is over all vectors $\alpha $ of nonnegative integers, and similarly
\begin{equation}\label{h shift}
h_{\lambda }(z_{1},z_{2},\ldots ,y_{1},y_{2},\ldots )=\sum _{\alpha }
h_{\alpha }({\bf z}) h_{\lambda -\alpha }({\bf y}) 
\end{equation}
or
$$ h_{\lambda }(a,y_{1},y_{2},\ldots )=\sum _{\alpha } a^{|\alpha |} h_{\lambda
-\alpha }({\bf y})$$
where the sum is over all vectors $\alpha $ of nonnegative integers, and
finally 
$$ \mbox{pow}_{\lambda }({\bf y}\cup {\bf z}) = \sum _{\mu \subseteq \lambda }
\prod _{i} { \mbox{mult}_{i }(\lambda ) \choose \mbox{mult}_{i}(\mu)}
\mbox{pow}_{\mu }({\bf y}) \mbox{pow}_{\lambda /\mu }({\bf z}).$$

Equations (\ref{e shift}) and (\ref{h shift}) are typical of full sequences of
divided powers. Let $q_{\lambda }({\bf y})$ be a full sequence with 
$$ q_{\lambda }({\bf y}\cup {\bf z})=\sum _{\alpha } q_{\alpha }({\bf y})
q_{\lambda -\alpha } ({\bf z}).$$
Then $q_{\lambda }({\bf y})$ is a {\em full sequence of divided
powers}.\index{Divided Powers} Note that $q_{(n)}({\bf y})$ is a linear
sequence of divided powers. Conversely, note that if $q_{n}({\bf y})$ is a
linear sequence of divided powers and their product $q_{\lambda }({\bf
y})=\prod_{i}q_{\lambda _{i}}({\bf y})$ is a full sequence, then it is a full
sequence of divided powers. However, that is not necessarily the only extension
of $q_{n}({\bf y})$ to a full sequence, and the product may not be full. For
example, 
$$ \prod _{i\geq 1}\id{\lambda _{i}}=\id{|\lambda |}. $$

Similarly, if $p_{\lambda }({\bf y})$ is a full sequence obeying
\begin{equation}\label{fully}
 p_{\lambda }({\bf y}\cup {\bf z})=\sum _{\alpha } {\lambda \choose \alpha
}p_{\alpha }({\bf y}) p_{\lambda -\alpha } ({\bf z})
\end{equation}
where ${\lambda \choose
\alpha }=\frac{\lambda !}{\alpha !(\lambda -\alpha )!}$ and $\alpha !=\prod
_{i}\alpha _{i}!$, then $p_{\lambda }({\bf y})$ is a {\em full sequence of
binomial type}\index{Binomial Type}. Its linear subsequence $p_{(n)}{\bf (y)}$
is a linear sequence of binomial type. 

Clearly, if $p_{\lambda }({\bf y})=\lambda !q_{\lambda }({\bf y})$, then
$p_{\lambda }({\bf y})$ is a full sequence of binomial type if and only if
$q_{\lambda }({\bf y})$ is a full sequence of divided powers.

\ssect{Transfer Operators and Adjoints}

There is a unique map which sends any full sequence to another full
sequence. In particular, we are interested in maps which send one full sequence
of divided powers (or equivalently binomial type) onto another. Such maps we
call {\em transfer operators}.\index{Transfer Operators}

\begin{thm}
\begin{enumerate}
\item The image of an full sequence of divided powers (resp. binomial type)
under any transfer operator is another transfer operator.
\item A linear operator on symmetric functions is a transfer operator if and
only if it is a homogeneous coalgebra isomorphism.
\end{enumerate}
\end{thm}

\proof{Coalgebra maps preserve the symmetric shift $E^{a}$.}
We can define the adjoint of any linear map on symmetric functions or
shift-invariant operators in a straightforward fashion; we are most
interested in the adjoints of transfer operators (that is, coalgebra maps), and
of automorphisms of the algebra of shift-invariant operators.
\begin{defn}\label{adjoints}
Let $\theta $ be a linear map on the space of shift-invariant operators, and
let $\phi $ be a linear map on the space of symmetric functions. Then
$\adj{\theta }$ and $\adj{\phi }$ are defined by the relations
\begin{eqnarray*}
\epsilon (\theta f({\bf D}))p({\bf y}) &=& \epsilon f({\bf D})(\adj{\theta
}p({\bf y}))\\* 
\epsilon (\adj{\phi }f({\bf D}))p({\bf y}) &=& \epsilon f({\bf D}) (\phi p({\bf
y}))
 \end{eqnarray*}
for all shift-invariant operators $f({\bf D})$ and symmetric functions $p({\bf
y})$. 
\end{defn}

By Theorems~\ref{sex} and \ref{stall}, these adjoints are well defined. For
example, the adjoint of the operator $f({\bf D} )$ on symmetric functions is
the operator which multiplies shift-invariant operators by $f({\bf D})$. 

These adjoints are adjoints in the usual sense
with respect to the action on the Hopf algebra of symmetric functions
of the Hopf algebra of shift-invariant operators which has been identified with
the dual of the Hopf algebra of symmetric functions.

Thus, $\adj{\adj{\theta }}=\theta $ and $\adj{\adj{\phi }}=\phi $. Moreover
the adjoint of any transfer operator is an automorphism of the algebra of
symmetric functions and visa versa. 
\begin{prop}\label{explitive}
The adjoint is an automorphism between the group of transfer operators,
and the group of automorphisms of the algebra of shift-invariant operators.\Box
\end{prop}

We can make this connection more explicit.
\begin{thm}\label{youth}
Let $\theta $ be a transfer operator which maps the linear sequence of binomial
type $p_{n}({\bf y})$ 
associated with the quasi-species $(a_{n})$ to the  linear sequence of binomial
type $q_{n}({\bf y})$ 
associated with the quasi-species $(b_{n})$.
Let $\adj{\theta }D_{n}/n!=\sum _{\lambda \vdash n}c_{\lambda }D_{\lambda }$.
Then 
$$ b_{n}= \sum _{\lambda \vdash n}c_{\lambda }\prod _{i\geq 1}
a_{\lambda_{i}}.$$   
\end{thm}

\proo{Consider the following sequence of equalities.
\begin{eqnarray*}
b_{n}&=& \epsilon \D _{n}q_{n}({\bf y})/n! \\
&=& \epsilon \D _{n}\theta p_{n}({\bf y})/n! \\
&=& \epsilon \adj{\theta } \D _{n}p_{n}({\bf y})/n!\\
&=& \sum _{\lambda \vdash n} c_{\lambda } \epsilon \D _{\lambda } p_{n}({\bf y})\\
&=& \sum _{\lambda \vdash n} c_{\lambda } \prod _{i\geq 1} a_{\lambda _{i}}\Box
\end{eqnarray*}}

\begin{thm}\label{Above}
Let $\theta $ be the linear operator on shift-invariant operators given by
$$ \theta \D _{\lambda}=\sum _{\mu \in {\cal P}} c_{\lambda \mu }\D _{\mu }, $$
and let $\phi $ be the linear operator on symmetric functions given by
$$ \phi m_{\lambda }({\bf y})= \sum _{\mu \in {\cal P}}d_{\lambda \mu } m_{\mu
}({\bf y}).$$ 
Then $\theta $ and $\phi $ are adjoint if and only if the matrices
$(c_{\lambda \mu })_{\lambda ,\mu \in {\cal P}}$ and $(d_{\lambda \mu
})_{\lambda ,\mu \in {\cal P}}$ are transposes of each other. That is,
$\adj{\theta }=\phi $ if and only if $c_{\lambda \mu }=d_{\mu \lambda }$ for
all partitions $\lambda $ and $\mu $.
\end{thm}

\proo{Consider the following series of equalities.
\begin{eqnarray*}
c_{\lambda \mu } &=& \epsilon \sum _{\nu \in {\cal P}} c_{\lambda \nu }\D_{\nu
} m_{\mu }({\bf y})\\ 
&=& \epsilon \adj{\phi } \D _{\lambda }m_{\mu }({\bf y})\\
&=& \epsilon \D _{\lambda } \phi m_{\mu }({\bf y})\\
&=& \epsilon \D _{\lambda} \sum _{\nu\in  {\cal P}}  d_{\mu \nu } m_{\nu }({\bf
y})\\ 
&=& d_{\mu \lambda }.\Box 
\end{eqnarray*}}

The classical involution $\omega  $ of the symmetric functions is the transfer
operator which sends the elementary symmetric function to the complete
symmetric function and back again. By Theorem~\ref{Above}, the coefficients of
$\adj{\omega  }$ are given by the coefficients of the {\em forgotten symmetric
function}\index{Forgotten Symmetric Function} $f_{\lambda }({\bf y})=\omega 
m_{\lambda }({\bf y})$. Thus, by \cite[Theorem 8(ii)]{Doubilet},
\begin{cor}\label{Doubilet}
The adjoint of the classical involution $\omega  $ of the symmetric functions
is given by the involution
$$ \adj{\omega  } \D _{n} = \sum _{\lambda\vdash n} \frac{n!}{\prod _{i}
i^{\lambda _{i}}} \D_{\lambda }.\Box $$
\end{cor}

\ssect{Genera}

We have left unanswered the question of what these full sequences of binomial
type enumerate. Surely, they do not count functions enriched with mere species,
for they do not possess enough structure.
Indeed, we define a generalization of the species called the {\em genus}
(plural: {\em genera}) whose enriched functions are counted by full sequences.

\begin{defn}\label{genus}\bold{Genera}
A genus is a functor {\bf G} from the category {\bf Part} of sets with partitions to
the category 
of sets. For any partition $\pi $ of a finite set $E$, we say that the members
of the set $[\pi] $ are 
${\bf G}$-structures, and for any bijection $f: E\rightarrow F$, we describe the
function  ${\bf G}[f]:{\bf G}[\pi ]\rightarrow {\bf G}[f(\pi) ]$ as a
relabeling of ${\bf G}$-structures. 
For this paper, we need to 
assume that there is only one structure on the empty partition and there is at
least one structure on partitions of a one element set; that is,
$|{\bf G}[\emptyset]|=1$ and $|{\bf G}[\{\{0 \} \}]|\neq 0$.
\end{defn}
Note that this definition is equivalent to the partitionals of
\cite{Nava}; however, it differs in the
multiplication, exponentiation,  and composition to be defined later.
Species are then be seen as the restriction of genera to set partitions
consisting of only one block.

The equivalence class of ${\bf G}$ is thus determined by the number $G_{\lambda
}$ of ${\bf G}$-structures on set partitions $\pi \vdash S$ of type $\lambda$. 
(The type of  $\pi $ is the integer partition $\lambda $ which
corresponds to the size of its blocks. For example, if $\pi =\{\{a,b,c \},\{d,e
\},\{f \},\{g \} \}$, then the type of $\pi $ would be $\lambda = (3211)$.)
The generating function of ${\bf G}$ is defined to be the formal power
series 
$$ \mbox{Gen}[{\bf G}]({\bf y}) = \sum _{\lambda \in {\cal P}} a_{\lambda
}m_{\lambda }({\bf y})/\lambda !.$$ 

As with species, again we consider the category  {\bf C}  of complex numbers.
\figx{{Quasi-Genera}\label{q genus}}
{$ \begin{array}{ccccc}
{\bf Sets} & {\iota \atop \displaystyle \longrightarrow} &{\bf Part}& {\SB{G}\atop
\displaystyle \longrightarrow} & {\bf Sets}\\* 
&&\thinr{\SB{Q}}\downarrow & \swarrow \thinl{\# }\\*
&&{\bf C}
\end{array} $}
Any functor from {\bf Part} to {\bf C} is called a {\em quasi-genus}. 
Note that any genus {\bf G} can be extended a quasi-genus {\bf Q} by composing
with $\#  $. We will perform computations with quasi-genera as if they were
actual genera. In 
this case, $a_{\lambda }$ above will refer to the value of $Q$ on a partition
of type $\lambda $.

Conversely, any genus can be restricted to a species by composing on the other
side with the functor $\iota $ from the category {\bf Sets} to the category
{\bf Part} defined by $\iota : E\mapsto \{E \}$.

The generating functions for a genus ${\bf G}$ and for its associated species
${\bf S}={\bf G}\circ\iota $ are simply related.
$$ \mbox{Gen}[{\bf G}](t,0,0,\ldots )=\mbox{Gen}[{\bf S}](t). $$
(See Example~\ref{nifty} for a converse.)

We define the sum of two genera ${\bf G}_{1}$ and ${\bf G}_{2}$ on a partition $\pi
\vdash S$ by their disjoint union
$$ ({\bf G}_{1}+{\bf G}_{2})[\pi ] = {\bf G}_{1}[\pi] \dot{\cup } {\bf G}_{2}[\pi ]$$
so that $$\mbox{Gen}[{\bf G}_{1}]({\bf y})+\mbox{Gen}[{\bf
G}_{2}]({\bf y}) = \mbox{Gen}[{\bf G}_{1}+{\bf G}_{2}]({\bf y}).$$
In other words, a ${\bf G}_{1}+{\bf G}_{2}$-structure is either a ${\bf
G}_{1}$-structure 
 or a ${\bf G}_{2}$-structure. With respect to this addition, the set of all
quasi-genera is a complex vector space. Hence, we can consider the following
pseudo-basis. 
$$ {\bf G}^{(\lambda )}[\pi ] = \left\{ \begin{array}{ll}
1&\mbox{if $\pi $ is of type $\lambda $, and}\\
0&\mbox{otherwise.}
\end{array} \right. $$
The additive identity then is the genus with no structures whatsoever on any
partition.

Similarly, we define the product of ${\bf G}_{1}$ and ${\bf G}_{2}$ on $\pi
\vdash S$ to be set of quadruples 
$$ ({\bf G}_{1}{\bf G}_{2})[\pi ]=
 \left\{ 
\parbox{3.5in}{$(S_{1},S_{2} ,A_{1},A_{2})$: $S_{1}$ and $S_{2}$ are disjoint,
$S_{1}\cup 
S_{2}=S,$ and $A_{i}\in {\bf G}_{i}[\pi _{i}]$ where $\pi _{i}$ is the
restriction of $\pi $ to $S_{i}$.}\right\}$$
In other words, as in Figure~\ref{Product G} in which each of the horizontal
lines represents a block of $\pi \vdash S$,  we divide $S$ in half---right and
left, \figx{Product of Genera\label{Product G}}{\vspace{2in}}
 and place a ${\bf G}_{1}$-structure on the
restriction of $\pi $ to one half, and a ${\bf G}_{2}$-structure on the restriction
of $\pi $ to the other half. 
\begin{prop}\label{times}
Let ${\bf G}_{1}$ and ${\bf G}_{2}$ be any two genera, then
$\mbox{Gen}[{\bf G}_{1}]({\bf y})\mbox{Gen}[{\bf G}_{2}]({\bf y}) =
\mbox{Gen}[{\bf G}_{1}{\bf G}_{2}]({\bf y})$. 
\end{prop}

{\em Proof:} By linearity, it will suffice to consider the case ${\bf
G}_{1}={\bf G}^{(\lambda )}$ and ${\bf G}_{2}={\bf G}^{(\mu )}$. Now, let
$\pi =\{\pi _{1},\pi _{2},\ldots  \}\vdash S$ be of type $\nu $ where the
blocks $\pi _{i}$ are listed in decreasing order of size. We must count
the number of ways to divide 
split $S$ into $S_{1}$ and $S_{2}$ such that $\pi $ restricted to $S_{1}$ and
$S_{2}$ is of type $\lambda $ and $\mu $ respectively.
Suppose further that $S_{1}\cap \pi _{i}$ contains $\alpha _{i}$ elements and
$S_{2}\cap \pi _{i}$ contains $\beta _{i}$ elements. Then $\alpha $ and $\beta
$ must be permutations of $\lambda $ and $\mu $ respectively, and $\alpha
+\beta =\nu $. Subject to these 
conditions then there are $\nu !/\lambda ! \mu !$ choices for $S_{1}$ and
$S_{2}$.
Hence,
$$ {\bf G}^{(\lambda )}{\bf G}^{(\mu )}= \sum _{\nu \in {\cal P}} \sum
_{\scriptstyle \alpha +\beta =\nu \atop {\scriptstyle \alpha \mbox{ is a
permutation of }\lambda \atop \scriptstyle \beta \mbox{ is a
permutation of }\mu }} \frac{\nu !}{\lambda !\mu !} {\bf G}^{(\nu )}.$$
On the other hand,
\begin{eqnarray*}
\frac{m_{\lambda }({\bf y})}{\lambda !} \frac{m_{\mu }({\bf y})}{\mu !}
&=& \sum _{\scriptstyle \alpha \mbox{ is a
permutation of }\lambda \atop \scriptstyle \beta \mbox{ is a
permutation of }\mu } {\bf y}^{\alpha +\beta }/\lambda !\mu !\\
&=& \sum _{\nu \in {\cal P}} \sum _{\scriptstyle \alpha +\beta =\nu \atop
{\scriptstyle \alpha \mbox{ is a 
permutation of }\lambda \atop \scriptstyle \beta \mbox{ is a
permutation of }\mu }} \frac{\nu !}{\lambda !\mu !} \frac{m_{\nu}({\bf y})}{\nu
!}. \Box 
\end{eqnarray*}

Thus, the multiplicative identity is the genus ${\bf G}^{(0)}$ which has one
structure on the unique partition of the empty set, but has no structures for
any other partitions.

Next, we define the exponentiation of a genus ${\bf G}$ on $\pi \vdash S$ to be
$$ \exp({\bf G})[\pi ]=
 \left\{ 
\parbox{3.5in}{$(\phi , (A_{B})_{B\in \phi })$: $\phi$ is a partition of $S$,
$(A_{B})_{B\in \theta }$ is a sequence indexed by the blocks of $\phi $
for each block $B$ of the partition $\phi $,  
$A_{B}\in {\bf G}[\pi _{B}]$ where $\pi _{B}$ is the restriction of
$\pi $ to $B$.}\right\}$$
That is, we divide $S$ into a number of parts and place a ${\bf G}$ structure on the
restriction of $\pi $ to each part. 
\figx{Exponentiation and Composition of Genera\label{Exponent G}}{\vspace{2in}}
In Figure~\ref{Exponent G}, the horizontal lines represent the blocks of $\pi $
and the circled areas represent the blocks $\phi $.
Since $\exp ({\bf G})$ is equivalent to $\sum _{n\geq 0}{\bf G}^{n}/n!$, we
immediately have 
$$\exp \left(\mbox{Gen}[{\bf G}]({\bf y}) \right)=\mbox{Gen}[\exp ({\bf
G})]({\bf y}).$$ 

Similarly, we define the composition of ${\bf G}_{1}$ with ${\bf G}_{2}$ on
$\pi\vdash S $ to be 
$$ \exp({\bf G})[\pi ]=
 \left\{ 
\parbox{3.5in}{$(\phi , (A_{B})_{B\in \phi },C)$: $\phi$ is a partition of $S$,
$(A_{B})_{B\in \phi}$ is a sequence indexed by the blocks of $\phi $
for each block $B$ of the partition $\phi $,  
$A_{B}\in {\bf G}[\pi _{B}]$ where $\pi _{B}$ is the restriction of
$\pi $ to $B$, and $C\in {\bf G}_{1}[\zeta ]$ where 
$\zeta $ is the partition of the blocks of $\phi$ defined by the equivalence
relation $B_{1}\sim B_{2}$ if and only if for every $P\in \pi $, we have
$|B_{1}\cap P|= |B_{2}\cap P|$.}\right\}$$
That is, we divide $S$ into a number of parts, place a ${\bf G}_{2}$-structure
on the restriction of $\pi $ to each parts. Finally we consider the partition
of the parts of $S$ classified according to their intersections with the blocks
of $\pi $  (as shown in
in Figure~\ref{Exponent G} by wavy lines), and we place a ${\bf
G}_{1}$-structure on it. 

Recall that the plethystic composition $p\circ q({\bf y}) $of two symmetric
functions $p({\bf y})$ 
and $q({\bf y})$ is defined as follows. We express $g({\bf y})$ as a sum of
monomials $g({\bf y})=\sum _{\alpha } c_{\alpha }{\bf y}^{\alpha }$, and
introduce a new set of variables ${\bf z}$ defined by $\prod_{i\geq
1}(1+z_{i}t) = \prod_{\alpha }(1+x^{\alpha }t)^{c_{\alpha }}$.  Finally, we
define $p\circ q({\bf y})=p({\bf z})$.

The generating function of a composition of genera can be expressed in term of
plethystic composition. 
\begin{thm}
Let ${\bf G}_{1}$ and ${\bf G}_{2}$ be any two genera, then
$\mbox{Gen}[{\bf G}_{1}({\bf G}_{2})]({\bf y}) = (\mbox{Gen}[{\bf G}_{1}]\circ
\mbox{Gen}[{\bf G}_{2}]) ({\bf y}).$ 
\end{thm}

\proof{By linearity, it will suffice to consider the case ${\bf G}_{1}={\bf
G}^{(\lambda )}.$ Now, let ${\bf G}= {\bf G}_{2}=\sum _{\mu \in {\cal P}}c_{\mu
} {\bf 
G}^{(\mu )}$. Now, let's compute the number of ${\bf G}^{(\lambda )}({\bf G})$
structures on the partition $\pi =\{\pi _{1},\pi _{2},\ldots  \}\vdash S$ of
type $\nu $ where
the blocks $\pi _{i}$ are listed in decreasing order of size. We must divide
$S$ into exactly $|\lambda |$ parts which intersect with $\pi $ in exactly
$\ell (\lambda )$ different ways with the $i\th$ way occuring exactly $\lambda
_{i}$ times. Let $\alpha ^{(1)},\alpha ^{(2)},\ldots
,\alpha ^{(\ell (\lambda ))}$ represent these various intersection types.
These intersection types must all be distinct. However, two sets of 
intersection types differing only by the rearrangement of $\alpha ^{(i)}$ and
$\alpha ^{(j)}$ corresponding to $\lambda _{i}=\lambda _{j}$ are not considered
distinct; this caveat will be denoted (*). Next,
let $S_{i1},S_{i2},\ldots ,S_{i\lambda _{i}}$ be the $\lambda _{i}$ blocks of
intersection type $\alpha ^{(i)}$. That is, $S_{ij}\cap \pi _{k}=\alpha
^{(i)}_{k}$ and the $S_{ij}$ form a partition of $S$. Then $\nu =\lambda
_{1}\alpha^{(1)}+\cdots +\lambda _{\ell(\lambda )}\alpha ^{(\ell (\lambda ))}$.
Each $S_{ij}$ is then endowed with a ${\bf G}$ structure in one of $c_{\mu }$
way where $\mu^{(i)} $ is the partition formed by sorting the vector $\alpha
^{(i)}$. Hence, 
\begin{eqnarray*}
{\bf G}^{(\lambda )}({\bf G}) 
&=&  \sum _{\nu}  \sum _{\mu ^{(i)}\in {\cal P}}
\sum _{\scriptstyle \lambda _{1}\alpha ^{(1)}+\lambda_{2}\alpha^{(2)}+\cdots
=\nu  \atop \scriptstyle (*), \alpha ^{(i)}\mbox{ is a permutation of }\mu
^{(i)}} 
\prod _{i=1}^{\ell (\lambda )}\left(\frac{c_{\mu ^{(i)}}}{\mu ^{(i)}!}
\right)^{\lambda _{i}} {\bf G}^{(\nu)} \nu!/ \lambda ! \\  
&=&  \theta m_{\lambda }({\bf H}^{(1)}({\bf y}),{\bf H}^{(2)}({\bf y}),\ldots )/\lambda !
\end{eqnarray*}
where the ${\bf H}^{(1)}({\bf y})$ are a set of genera of which $c_{\mu }/\mu!$
equals ${\bf G}^{(\mu)}$ for each partition $\mu \in {\cal P}$ and $\theta $ is
the 
linear operator $\theta :m_{\nu }({\bf y})\mapsto \nu !m_{\nu }({\bf y})$.  The
conclusion now follows immediately.} 

Thus, the compositional identity is the genus ${\bf G}^{(1)}$ which has one
structure on the unique partition of any one element set, and no structures
otherwise. 

Finally, we define the $n\th $ iterated derivative ${\bf D}_{n}{\bf G}$of a
genus ${\bf G}$ on a partition $\pi$ to be the set consisting structures of
${\bf G}$-structures on a partition $\pi \cup \{\{\infty _{1},\ldots ,\infty
_{n} \} \}$ containing a new block of size $n$ along with a complete ordering
of this new block.
Thus, ${\bf D}_{n}{\bf G}^{(\lambda) }=n!{\bf G}^{(\lambda \backslash n)}$ so
that 
$$\D_{n}\mbox{Gen}[{\bf G}]({\bf y}) = \mbox{Gen}[{\bf D}_{n}{\bf G}]({\bf
y}).$$

By considering species and their generating functions, we immediately observe
that ${\rm pow}_{n}\circ{\rm pow}_{m}({\bf y})={\rm pow}_{nm}({\bf y})$ from
which we can deduce that general composition is associative.

Having now completely defined an adequate generalization of species. We can now
consider function enriched by it.
\begin{defn}[Generic Functions]\label{Generic}
Let $S$ and $X$ be sets, and let $\pi $ be a partition of $S$. $f$ is called a
generic 
function or a function enriched by the genus ${\bf G}$ if it is a function from $S$
to $X$ which is equipped with a ${\bf G}$-structure for the restriction of $\pi $ to
each of its fibers.
\end{defn}
\figx{Typical Generic Function}{\vspace{2in}}

\begin{thm}
A full sequence of symmetric functions $p_{\lambda }({\bf y})$ is of binomial
type if and only if it enumerates the generic functions from a set $S$ equipped
with a partition $\pi \vdash S$ of type $\lambda $ to the set
$X=\{y_{1},y_{2},\ldots  \}$ for some genus (or quasi-genus) ${\bf G}$.
\end{thm}

\proof{{\bf (If)} The left side of \eref{fully} enumerates the set of generic
functions from $\pi \vdash S$ to $X_{1}\cup X_{2}$. The right side of
\eref{fully} counts the number of ways to divide $S$ into two parts, and map
the parts to $X_{1}$ and $X_{2}$ respectively with generic functions. The sum 
is over the types of the resulting restrictions of the partition $\pi $.
Obviously, these two set of generic functions are identical.

{\bf (Only If)} Define the numbers ${\bf G}_{\lambda ,\alpha }$ by the identity
$$ D_{n}p_{\lambda }({\bf y}) = \sum _{\alpha } {\lambda \choose \alpha }
{\bf G}_{\lambda ,\alpha }p_{\lambda -\alpha }({\bf y}) $$.
Thus,
$$ E^{z}p_{\lambda }({\bf y}) = \sum _{\alpha } {\lambda \choose \alpha  }
\frac{z^{|\alpha |{\bf G}_{\lambda ,\alpha }}}{|\alpha |!} p_{\lambda -\alpha ({\bf
y})}. $$
However, by \eref{fully}, ${z^{|\alpha |{\bf G}_{\lambda ,\alpha }}}/{|\alpha |!}
$ does not depend on $\lambda $. Hence, ${\bf G}_{\lambda ,\alpha }={\bf G}_{\alpha }$ does
not depend on $\lambda $. Moreover, it is invariant under permutation of
$\alpha $. Thus, ${\bf G}_{\alpha }$ is a quasi-genus, and $p_{\lambda }({\bf y})$
enumerates its generic functions.}

Restricting our attention to species, we have
\begin{cor}
A linear sequence of symmetric functions $p_{n }({\bf y})$ is of binomial
type if and only if it enumerates the functions from $N$ to
$X=\{y_{1},y_{2},\ldots  \}$ enriched by  some species (or quasi-species) ${\bf
S}.\Box$ 
\end{cor}

The two examples of full sequences of binomial type given are associated with
fairly simple genera.
\begin{ex}
$e_{\lambda }({\bf y})$ enumerates functions enriched with a genus ${\bf G}$
with a 
single structure for any partition consisting solely of one element blocks, and
no structures otherwise.
\end{ex}
\begin{ex}
$h_{\lambda }({\bf y})$ enumerates functions enriched with a genus ${\bf G}$ with
$\lambda !$ structures on any partition of type $\lambda $. For example,
${\bf G}[\pi ]$ could be the set of posets which completely order the blocks of
$\pi $ yet leave elements of distinct blocks incomparable.
\end{ex}

As we have explained about the products of most linear sequences of symmetric
functions of binomial type form similar full sequences.
\begin{ex}\label{nifty}
Let $p_{n}({\bf y})$ be the linear sequence of symmetric functions of binomial
type associated with the quasi-species ${\bf S}=(a_{n})_{n\geq 0}$. Then if
full, $p_{\lambda 
}({\bf y}) =\prod _{i}p_{\lambda _{i}}({\bf y})$ is the full sequence of
symmetric functions of binomial type associated with the quasi-genus
${\bf G}=(a_{\lambda })_{\lambda \in {\cal P}}$ where $a_{\lambda }=\prod_{i}
a_{\lambda _{i}}$. That is, 
each ${\bf G}$ structure on $\pi $ is a sequence of $S$-structures---one for each
block of $\pi $. The generating function for ${\bf G}$ is
$$ \mbox{Gen}[{\bf G}]({\bf y})=\prod _{y\in X} \mbox{Gen}[{\bf S}](y). $$
\end{ex}

Conversely, we have noted that every full sequence of
symmetric function of binomial type restricts to a similar linear sequence
\begin{ex}
Let $p_{\lambda }({\bf y})$ be the full sequence of symmetric functions of
binomial type associated with the quasi-genus ${\bf G}=(G_{\lambda })_{\lambda \in
{\cal P}}$. Then $p_{n}({\bf y})=p_{(n)}({\bf y})$ is the linear sequence
associated with the quasi-species ${\bf S}={\bf G}\circ\iota $. That is, $S_{n}=S_{(n)}$,
and for any set $E$, $S[E]=G[\{E \}].$ (See Figure~\ref{q genus}.)
\end{ex}

Thus, derivatives of a full sequence of symmetric functions of binomial type
are now easy to compute.
\begin{prop}\label{duffer}
Let $p_{\lambda }({\bf y})$ be the full sequence of symmetric functions of
binomial type associated with the quasi-genus $(G_{\lambda })_{\lambda \in
{\cal P}}$. Its symmetric derivative is
$$ \D p_{\lambda }({\bf y}) = \sum _{i\geq 1}iG_{(1)}\mbox{mult}_{i}(\lambda )
p_{\mu }({\bf y})$$
where $\mu =\lambda \backslash i \cup \{i-1\}$,
and more generally
$$ \D _{n} p_{\lambda }({\bf y})= n! \sum _{|\alpha |=n} {\lambda \choose \alpha
}G_{\alpha }p_{\lambda-\alpha }$$
\end{prop}

\proof{To compute the $n\th$ iterated derivative, we distinguish an element of
the range of the generic functions,\figx{Contribution to the Derivative of a
Full Sequence of Binomial Type}{\vspace{2in}} and count generic functions for
which the inverse image of that 
element contains exactly $n$ points. Here, we sum over how these $n$ points are
arranged within the block structure of the domain of the generic function.}

\begin{cor}\label{shifts}
Let $p_{\lambda }({\bf y})$ be the full sequence of symmetric functions of
binomial type associated with the quasi-genus ${\bf G}=(G_{\lambda })_{\lambda
\in {\cal P}}$. Its shift is 
$$ E^{x}p_{\lambda }({\bf y}) = \sum _{\alpha } { \lambda \choose \alpha  }
x^{|\alpha |}G_{\alpha }p_{\lambda-\alpha }({\bf y}).$$
\end{cor}

From Corollary~\ref{shifts} and \eref{fully}, we have the following explicit
formula for a full sequence of symmetric functions of binomial type.
\begin{cor}\label{bliss}
The full sequence of symmetric functions of
binomial type $p_{\lambda }({\bf y})$ associated with the quasi-genus
$(G_{\lambda })_{\lambda \in {\cal P}}$ is given by the formula
$$ p_{\lambda }(x,y,z,\ldots )= \sum _{\lambda =\alpha +\beta +\gamma +\cdots }
G_{\alpha }G_{\beta }G_{\gamma }\cdots x^{|\alpha |}y^{|\beta |}z^{|\gamma
|}\cdots  .\Box $$
\end{cor}
Thus, its evaluation is given by
$$ \epsilon p_{\lambda }({\bf y})=\delta _{\lambda ,(0)}. $$
Also, we see that the definition of a full sequence was not really as strict as
one may have thought.
\begin{cor}
If $p_{\lambda }({\bf y})$ is associated with a quasi-genus, and is a basis for
$\Lambda $, then it is a full sequence.$\Box$
\end{cor}

Now, by iterating Proposition~\ref{duffer}, we have
\begin{cor}\label{repeat}
Let $p_{\lambda }({\bf y})$ be the full sequence of symmetric functions of
binomial type associated with the quasi-genus $(G_{\lambda })_{\lambda \in
{\cal P}}$. Then
$$ \D _{\mu }p_{\lambda }({\bf y}) = \sum _{|\alpha ^{(i)}|=\mu _{i}} {\lambda
\choose \alpha^{(1)},\alpha ^{(2)},\ldots \in {\cal P} } G_{\alpha ^{(1)}}G_{\alpha
^{(2)}}\cdots  p_{\lambda -\alpha ^{(1)}-\alpha
^{(2)}} ({\bf y})$$ 
where each $\alpha ^{(i)}$ is a partition, and ${\lambda \choose \alpha
^{(1)},\alpha ^{(2)},\ldots  }$ is given by the product of multinomial
coefficients
$${\lambda \choose \alpha
^{(1)},\alpha ^{(2)},\ldots  }=\prod_{i\geq 1} {\lambda _{i}\choose \lambda
_{i}-\alpha _{i}^{(1)}-\alpha _{i}^{(2)}-\cdots ,\alpha
_{i}^{(1)},\alpha _{i}^{(2)},\ldots  } =\frac{\lambda !}{(\lambda -\alpha
^{(1)}-\alpha ^{(2)}-\cdots )!\alpha ^{(1)}!\alpha ^{(2)}!\cdots } .\Box$$
\end{cor}

Finally, we compute the coefficients of transfer operators in terms of genera.
\begin{thm}\label{age}
Let $p_{\lambda }({\bf y})$  be the full sequence
of binomials type associated with the quasi-genus ${\bf G}=(G_{\lambda
})_{\lambda \in {\cal P}}$. Suppose $\theta $ is the transfer operator given by
$$ \adj{\theta } \D _{n}/n! = \sum _{\lambda \vdash n}c_{\lambda }\D_{\lambda
}.$$ 
Then $\theta p_{\lambda }({\bf y})$ is the full sequence of binomial type
associated with the quasi-genus $(H_{\lambda })_{\lambda \in {\cal P}}$ where
$$ H_{\lambda }=\sum _{\mu \vdash n}c_{\mu }\sum _{M} {\lambda \choose
M} G_{M_{1}}G_{M_{2}}\cdots  $$
where $M$ is a matrix of nonnegative integers whose column sums are $\lambda $,
whose row sums are $\mu $, and whose $i\th$ row is $M_{i}$ and ${\lambda
\choose M}=\lambda !/\prod _{i,j\geq 1}M_{ij}!$.
\end{thm}

\proo{Consider the following sequence of equalities.
\begin{eqnarray*}
H_{\lambda }&=& \epsilon \sum _{|\alpha |=n} {\lambda \choose \alpha }H_{\alpha
}q_{\lambda -\alpha }({\bf y})\\
&=& \epsilon \D_{n} q_{\lambda }({\bf y})/n! \\
&=& \epsilon \D_{n} \theta p_{\lambda }({\bf y})/n! \\
&=& \epsilon \adj{\theta} \D _{n} p_{\lambda }({\bf y})/n! \\
&=& \epsilon \sum _{\mu \vdash n} c_{\mu } \D_{\mu } p_{\lambda } ({\bf y})\\
&=& \epsilon \sum _{\mu \vdash n} c_{\mu } \sum _{|\alpha ^{(i)}|=\mu _{i}}
{\lambda \choose \alpha ^{(1)},\alpha ^{(2)},\ldots  } G_{\alpha
^{(1)}}G_{\alpha ^{(2)}}\cdots p_{\lambda -\alpha ^{(1)-\alpha ^{(2)}-\cdots }}({\bf y})\\
&=& \sum _{\mu \vdash n}c_{\mu }\sum _{M} {\lambda \choose
M} G_{M_{1}}G_{M_{2}}\cdots  \Box 
\end{eqnarray*} }

When one takes $\lambda =(n)$, Theorem~\ref{age} specializes to
Theorem~\ref{youth}.

\end{document}